\makeatletter
\IfFileExists{./numapde-latex/numapde-packages.sty}{\providecommand*{\input@path}{}\edef\input@path{{./numapde-latex/}\input@path}}{}
\makeatother

\documentclass{numapde-preprint}

\usepackage{numapde-semantic}
\usepackage{numapde-style}
\usepackage{numapde-local}

\IfFileExists{./numapde-bibliography/numapde.bib}{\addbibresource{./numapde-bibliography/numapde.bib}}{\addbibresource{numapde.bib}}

\hypersetup{
	pdftitle={Numerical investigation of agent controlled pedestrian dynamics using a structure preserving finite volume scheme},
	pdfauthor={Jan-Frederik Pietschmann, Ailyn Stötzner, Max Winkler},
	pdfkeywords={crowd motion, nonlinear transport, Eikonal equation, ODE-PDE coupling, optimal control, finite volume, projected gradient descent}
}

\title{Numerical investigation of agent controlled pedestrian dynamics using a structure preserving finite volume scheme}
\subtitle{}
\shorttitle{Pedestrian dynamics}

\author{Jan-Frederik Pietschmann\thanks{Technische Universität Chemnitz, Faculty of Mathematics, AG Inverse Probleme, 09107 Chemnitz, Germany (\email{jan.pietschmann@mathematik.tu-chemnitz.de}, \url{https://www.tu-chemnitz.de/mathematik/invpde/prof.php}, \orcid{0000-0003-0383-8696}).}
\and
Ailyn Stötzner\thanks{Technische Universität Chemnitz, Faculty of Mathematics, 09107 Chemnitz, Germany (\email{ailyn@stoetzner.de}, \url{http://www.tu-chemnitz.de/mathematik/part_dgl/people/stoetzner/}, \orcid{0000-0002-8287-9738}, \email{max.winkler@mathematik.tu-chemnitz.de}, \url{https://www.tu-chemnitz.de/mathematik/part_dgl/people/winkler/}, \orcid{0000-0002-5292-2280}).}
\and
Max Winkler\footnotemark[2]}
\shortauthor{J. Pietschmann, A. Stötzner and M. Winkler}

\dedication{}

\IfFileExists{numapde-uncolorizeHyperrefIfChanges.sty}{\RequirePackage{numapde-uncolorizeHyperrefIfChanges}}{}

\begin{document}
\maketitle

\begin{abstract}
We provide a numerical realisation of an optimal control problem for pedestrian motion with agents that was analysed in \cite{2011.03580}. The model consists of a regularized variant of Hughes' model for pedestrian dynamics coupled to ordinary differential equations that describe the motion of agents which are able to influence the crowd via attractive forces. We devise a finite volume scheme that preserves the box constraints that are inherent in the model and discuss some of its properties. We apply our scheme to an objective functional tailored to the case of an evacuation scenario. Finally, numerical simulations for several practically relevant geometries are performed.\end{abstract}


\begin{AMS}
\href{https://mathscinet.ams.org/msc/msc2010.html?t=49K20}{49K20}, \href{https://mathscinet.ams.org/msc/msc2010.html?t=35Q91}{35Q91}, \href{https://mathscinet.ams.org/msc/msc2010.html?t=35M33}{35M33}, \href{https://mathscinet.ams.org/msc/msc2010.html?t=65M08}{65M08}
\end{AMS}

\section{Introduction}
\label{sec:introduction}

With more and more people living in highly populated areas, the modelling, simulation and control of (large) pedestrian crowds is an important field of research. In this work, we study the optimal control problem for a regularized version of Hughes' model for pedestrian motion, \cite{Hughes:2002:1}. In our approach, the (continuous) crowd can be controlled by a fixed small number of agents that can attract people in their vicinity. In terms of the model, this corresponds to an additional potential term centred at the agents positions. In a previous work, \cite{2011.03580}, we already studied the well-posedness and optimality conditions of this problem while here, we focus on a numerical implementation of the control problem and extensive numerical examples. In particular, we provide and analyse a finite volume scheme that preserves the box constraints inherent in our problem.

To introduce the model, we fix $\Omega\subset\R^2$ to be a bounded domain with
$C^4$-boundary $\partial\Omega$. Furthermore $T>0$ is an arbitrary time horizon and $Q_T:=(0,T)\times \Omega$ as well as $\Sigma_T = (0,T)\times \partial \Omega$ denote the
space-time cylinder and its lateral boundary, respectively. The boundary is decomposed
into two parts: $\partial\Omega_D$ representing the 
exits and $\partial\Omega_W$ the part where the domain is constrained by
walls. For theoretical purposes (regularity of solutions) we assume 
$\overline{\partial\Omega_D}\cap \overline{\partial\Omega_W}= \emptyset$
meaning that both boundary parts are separated from each other, see
\cref{fig:solution_example_1}. In a similar way we define $\Sigma_D = (0,T)\times
\partial\Omega_D$ and $\Sigma_W = (0,T)\times \partial\Omega_W$.\\
The unknown variables in our system of equations are the density of the crowd
$\rho\colon Q_T\to \R_+$, a potential specifying the current time to escape
$\phi\colon Q_T\to \R$.
In addition, there are $M$ agents which may influence the motion of the crowd via attractive forces. Their positions are denoted by $x_i\colon (0,T)\to
\R^2$, $i=1,\ldots,M$. In addition, each agent is able to regulate the strength by which it acts on the crowd. This is encoded in the intensities $c_i\colon (0,T) \to \R_+$, $i=1,\ldots,M$. Both the agent trajectories and interaction strength are summarized in a vector $\bx
= (x_1,\ldots,x_n)^\top$ and $\bc
= (c_1,\ldots,c_n)^\top$, respectively.

The mathematical equations describing the movement of a pedestrian crowd
influenced by agents then read as follows. For given agent movement directions
$\bu = (u_1,\ldots,u_M)^\top$ with $u_i\in L^\infty(0,T;\R^2)$ the unknowns
$\rho,\phi,\bx$ are related to each other by means of
\begin{subequations}
	\label{eq:forward_system}
	\begin{alignat}{2}
		\label{eq:forward_system1}
		\partial_t\rho - \nabla\cdot\paren[big](){\rho\,\beta(\rho,\phi,\bx,\bc)} 
		&
		= 
		\eps \, \laplace \rho
		&
		&
		\quad
		\text{in } Q_T,
		\\
		\label{eq:forward_system2}
		-\delta_1\,\laplace \phi + \abs{\nabla\phi}^2 
		&
		= 
		\frac{1}{f(\rho)^2+\delta_2}
		&
		&
		\quad
		\text{in } Q_T, 
		\\
		\label{eq:forward_system3}
		\dot x_i(t) 
		&
		=  
		f\paren[big](){\rho(t,x_i(t))} \, u_i(t)
		&
		&
		\quad
		\text{for } t\in(0,T), \quad i=1,\ldots,M.
	\end{alignat}
\end{subequations}
Moreover, we impose the boundary conditions
\begin{equation}\label{eq:forward_bc}
	\begin{aligned}
		-\paren[big](){\eps\,\nabla\rho  + \rho\,\beta(\rho,\phi,\bx,\bc)} \cdot n 
		&
		=
		\gamma\,\rho 
		,\quad
		& 
		\phi 
		&
		= 
		0 
		&
		& 
		\text{on } \Sigma_\textup{D}
		,
		\\
		\paren[big](){\eps\,\nabla\rho  + \rho\,\beta(\rho,\phi,\bx,\bc)} \cdot n 
		&
		=
		0 
		,
		& 
		\nabla \phi \cdot n 
		&
		=
		0 
		&
		&
		\text{on } \Sigma_\textup{W}
		,
	\end{aligned}
\end{equation}
as well as the initial conditions
\begin{equation}\label{eq:forward_initial_conditions}
	\rho(0,\cdot) = \rho_0\quad \text{in}\ \Omega,\quad x_i(0) = x_{i,0}\quad
	\text{for}\ i=1,\ldots,M.
\end{equation}
Here, $\varepsilon, \delta_1, \delta_2 > 0$ are regularization parameters and
the corresponding terms in the system are needed to guarantee a certain regularity for the
solution, see \cref{thm:ex_state}. 
The domain $\Omega$ is sufficiently large such that $x_i(t)\in \Omega$
on $[0,T]$ for $i=1,\ldots,M$ and $t\in [0,T]$ if $\abs{u_i(t)}\le 1$.

Let us briefly discuss the meaning of the respective terms: Equation
\eqref{eq:forward_system1} states that pedestrians are transported
according to the velocity field $\beta$, see \eqref{eq:def_beta}
below, while also performing (little) random motion encoded by the
Laplacian of $\rho$. The second equation \eqref{eq:forward_system2} is
a modified and regularized Eikonal equation whose solution is the
distance to the closest exit, mitigating areas of high density via
the term on the right-hand side. Here, the additional diffusion
accounts for the fact that pedestrians do not know their environment
exactly. Then, \eqref{eq:forward_system3} governs the motion of the
agents, whose speed is also influenced by the surrounding pedestrian
density.
The function $f\colon [0,1]\to [0,1]$ is
a density-velocity rule, chosen in such a way that $f(\rho)$
determines the maximum velocity an individual can move if the density
in its current position is $\rho$. We choose $f$ to be monotonically
decreasing meaning that higher densities lead to slower movements.
 The velocity field $\beta$ will reflect the fact that pedestrians are, on the one hand, trying the minimize their exit time which amounts to a drift term in the direction of $\nabla \phi$ and on the other hand, they are attracted by the agents which is realized by additional attractive potentials whose center depends on the agents' positions $\bx$. This results in a velocity which is the sum of two terms. Furthermore, to account for the effect that the velocity will deteriorate in regions of high density, it will be modified by an additional multiplicative factor $f(\rho)$. As in the equations for the motion of the agents, $f$ is monotonically decreasing and becomes zero at a given maximal density.

The boundary conditions \eqref{eq:forward_bc} allow for an outflow with velocity $\gamma$ on parts of the boundary ($\Sigma_D$) while no-flux conditions on the remaining parts are to be interpreted as walls ($\Sigma_W$). A detailed description of the involved non-linearities will be given in the next section, but we also refer to \cite{2011.03580} for more details on the model and the regularizing terms.


Analytical properties of the unregularized Hughes' model introduced in \cite{Hughes:2002:1} (i.e. $\eps = \delta_1 =\delta_2 =0$ in \eqref{eq:forward_system}), without control, are difficult because of the low regularity of $\nabla \phi$ that appears on a set depending on the solution $\rho$ of the first equation, but see \cite{AmadoriGoatinRosini:2014:1,AmadoriDiFrancesco:2012:1,ElKhatibGoatinRosini:2013:1}. Thus, regularized variants have been considered, see \cite{DiFrancescoMarkowichPietschmannWolfram:2011:1} for an instance where $\eps = 0$ but $\delta_1, \,\delta_2 \neq 0$. In fact, the result there is obtained as a vanishing viscosity limit $\eps \to 0$. There is also a number of extensions and variants of the model, aiming to understand additional properties, make it more realistic, or consider different settings like graphs, see \cite{BurgerDiFrancescoMarkowichWolfram:2014:1,CarrilloMartinWolfram:2016:1,CarliniFestaSilvaWolfram:2016:1,DiFrancescoFagioliRosiniRussoAnd:2017:1,ColomboGokieliRosini2018}.\\
Control of systems by means of a small number of agents has received lots of interest recently, both on a discrete level (i.e. one considers a large system of ODEs for the motion of individuals coupled to a small number of equations for the agents), see \cite{CaponigroFornasierPiccoliTrelat2013,HimakalasaWongkaew2021}, but also for coupled PDE-ODE systems, \cite{Albi2016,Albi2017}. Let us emphasize that whenever PDEs are coupled to ODEs in such a fashion that the solution of the PDE needs to evaluated at the solution of the ODE (as in \eqref{eq:forward_system3}), regularity is needed. While in our case, this is obtained by the additional diffusion in \eqref{eq:forward_system1}, when hyperbolic models for the transport of pedestrian are considered, an additional regularization in the ODE is needed, see e.g. \cite{Borsche:2015:1,BorscheKlarKuehnMeurer:2014:1,Borsche:2019:1}. Finally, let us mention that the ODE-ODE and PDE-ODE perspective are closely related by means of so-called mean field limits when the number of agents tends to infinity, see \cite{BurgerPinnauRothTotzeckTse:2016:1,PinnauTotzeck:2018:1} and also the recent overview on control of crowds, \cite{BandaHertyTrimborn2020}. \\
For the numerical discretization of \eqref{eq:forward_system1}, we employ, as
we think of the parameter $\eps$ being small, a finite volume scheme for the
spatial discretization, which may also be interpreted as a discontinuous
Galerkin scheme. In combination with the Lax-Friedrichs numerical fluxes the
scheme is stable and preserves the bounds $0 \le \rho \le 1$
inherent in our model. Such structure-preserving discretizations of PDEs
gained much attention, e.g., in the context of chemotaxis problems 
\cite{EpshteynKurganov2009,LiShuYang2017,StrehlSokolovKuzminTurek2010,GuoLiYang2018,IbrahimSaad2014,Filbet2006}.
The previously mentioned articles differ also in the choice of the
time-stepping scheme. For the treatment of \eqref{eq:forward_system} we will use an implicit-explicit finite
difference scheme, whereas the diffusion-related terms are established implicitly and
the convection-related terms explicitly. The equation
\eqref{eq:forward_system2} is discretized with standard linear finite elements
and for \eqref{eq:forward_system3} we employ a backward Euler scheme.

The paper is organized as follows: In \cref{sec:continuous} we give a precise definition of our problem and recall the analytical results from \cite{2011.03580}. In \cref{sec:num} we provide a numerical discretization scheme in space and time and analyse some of its properties, in particular, we show that it preserves physical bounds of the density of pedestrians. Corresponding optimization algorithms are discussed in \cref{sec:num_opt} and \cref{sec:num_ex} finally provides the results of our numerical experiments.

\section{The continuous optimal control problem}\label{sec:continuous}

Let us motivate the remaining quantities arising in the system of equations.
The function $f\colon[0,\rho_{\textup{max}}]\to \R$ is a density-velocity relation and is assumed
to be $W^{3,\infty}(\R)\cap C_c(\R)$ with $f(0)=1$ and
$f(\rho_{\textup{max}})=0$
with $\rho_{\textup{max}}$ denoting the maximal density. A usual
choice is
\begin{equation*}
	f(\rho) = \xi\,\left(1-\frac{\rho}{\rho_{\textup{max}}}\right)
\end{equation*}
with a cut-off function $\xi\in C_c^\infty(-1,2)$ satisfying $\xi\equiv 1$
on $(0,1)$. Obviously, a higher density leads to a lower velocity. Throughout the article we set
$\rho_{\textup{max}} =1$.
The movement direction of the crowd described by the function
$-\beta(\rho,\phi,\bx,\bc)$ is modelled as follows. The primary interest of the crowd is to move either towards the closest emergency exit, this is the
direction $-\nabla \phi$. This is mitigated by the attraction of close by agents which is the direction
$-\nabla \phi_K$, where is an agent potential defined by
\begin{equation*}
	\phi_K(\bx,\bc; t, x)
	\coloneqq 
	\sum_{i=1}^M c_i(t)\,K\paren[big](){x - x_i(t)}.
\end{equation*}
Here, $K(x) = k(\abs{x})$, $k\in W^{3,\infty}(\R)$ is a radially
symmetric function and $c_i\in H^1(0,T)$, $i=1,\ldots,M$, are
time-dependent intensity functions. Typical choices for attractive agent potentials are either the bumb function
\begin{equation*}
	k(r) = \begin{cases}
	     \exp\paren[auto](){-\frac{R^2}{R^2-r^2}}, &\text{if}\ r < 
	     R,\\
	     0, &\text{otherwise}
	\end{cases}
\end{equation*}
with attraction radius $R>0$, or the Morse potential
\begin{equation*}
	k(r) = \e^{-2\,a\,(r-r_a)} - 2\,\e^{-a\,(r-r_a)}
\end{equation*}
with certain parameters $a, r_a > 0$, realizing a repulsion in the near and an attraction in the far range of the agents. This is useful to avoid a high density very close to the respective agent. We refer to \cite{CARRILLOHUANGMARTIN2014} for a more detailed discussion on potentials in the context of flocking problems.

These considerations yield a velocity field defined as follows
\begin{equation}
\label{eq:def_beta}
	\beta(\rho,\phi,\bx,\bc) = v_0f(\rho)\,h(\nabla(\phi+\phi_K(\bx,\bc)))\quad
	\text{with}\quad h(x) = \min{}_{\varepsilon}\{1,\abs{x}\}\frac{x}{\abs{x}},
\end{equation}
where $h$ is a smoothed projection into the unit ball in $\R^2$ and the factor $f(\rho)$
again links the allowed movement speed to the current density, this is,
$\abs{\beta(\rho,\phi,\bx,\bc)} \le v_0f(\rho)$ in $Q_T$.

We briefly introduce the function spaces used in the sequel,
see also \cite{DenkHieberPruess2007}. For a
domain $\Omega\subset \R^2$ we denote by $W^{k,p}(\Omega)$, $k\in \N_0$, $p\in
[1,\infty]$, the usual Sobolev spaces and by $W^{k-1/p,p}(\Gamma)$ for $k\ge
1$ the corresponding trace spaces which may be equipped with the
Sobolev-Slobodetskij norm. Furthermore, we write $H^k(\Omega) =
W^{k,2}(\Omega)$. Special spaces incorporating already boundary
conditions are $H_D^1(\Omega) = \{v\in  H^1(\Omega) \colon
v|_{\partial\Omega_{\textup{D}}} \equiv 0\}$ and
$W^{2,p}_{\textup{DN}}(\Omega) := \{v\in W^{2,p}(\Omega)\colon
v|_{\partial\Omega_D} = 0,\ \partial_n v|_{\partial\Omega_W}=0\}$. 
For time-dependent functions $v\colon [0,T]\to X$
for some Banach space $X$ we define
\begin{equation*}
	L^p(0,T;X):=\{v\colon (0,T)\to X\ \vert\ \int_0^T \norm{v(t)}_X^p \d t < \infty\}, \quad p\in [1,\infty),
\end{equation*}
as well as
\begin{equation*}
	W^{s,p}(0,T;X):=\{v\colon (0,T)\to X\ \vert\ \partial_t^\ell
	v\in L^p(0,T;X), 0\le\ell\le s\},\quad s\in \N_0,\ p\in [1,\infty).
\end{equation*}
For the application we have in mind the following spaces
\begin{equation*}
	W^{r,s}_p(Q_T) 
	\coloneqq
	L^p(0,T;W^{r,p}(\Omega)) \cap W^{s,p}(0,T;L^p(\Omega))
	,\quad p \in [1,\infty),\quad r,s \in \N_0,
\end{equation*}
are of interest which are equipped with the natural norms
\begin{equation*}
	\norm{v}_{W^{r,s}_p(Q_T)} := \paren[auto](){\norm{v}_{L^p(0,T;W^{r,p}(\Omega))}^p +
	\norm{v}_{W^{s,p}(0,T;L^p(\Omega))}^p}^{1/p}.
\end{equation*}
Spaces with non-integral $r$ and $s$ are defined as interpolation spaces. 

In a previous work, \cite{2011.03580}, a global (in time) well-posedness and regularity result for
\eqref{eq:forward_system}--\eqref{eq:forward_initial_conditions} was established. Furthermore, optimality conditions for related optimal control problems where this system occurs as a constraint were derived. For convenience of the reader we briefly summarize the most important results needed in the present article.

First, there holds the following existence and regularity result:

\begin{theorem}\label{thm:ex_state} 
	Assume that $\rho_0\in W^{3/2,4}(\Omega)$ and fix $T>0$. Given arbitrary agent movement directions $\bu = (u_1,\ldots,
	u_M)^\transp \in L^\infty(0,T;\R^2)^M$ and intensities $\bc = (c_1,\ldots,c_M)\in H^1(0,T)^M$, there exists a
	unique strong solution $(\rho,\phi,\bx)$ to
	\eqref{eq:forward_system}--\eqref{eq:forward_initial_conditions}
	which satisfies, for any $2 <p < \infty$, $\rho \in W^{2,1}_p(Q_T)$ and $\phi\in
	L^\infty(0,T;W^{2,p}(\Omega))$. The agent trajectories $x_i$, $i=1,\ldots,M$ are
	absolutely continuous. 
	Moreover, the a~priori estimate
	\begin{equation*}
		\norm{\rho}_{W^{2,1}_p(Q_T)}
		+ \norm{\phi}_{L^\infty(0,T;W^{2,p}(\Omega))} 
		\le 
		C \norm{\rho_0}_{W^{1,p}(\Omega)},
	\end{equation*}
	holds with $C>0$ depending only on the domain, the bounds for the coefficients and the respective kernel.
\end{theorem}
The previous result allows to define an operator, the so-called
control--to--state operator,
\begin{equation*}
	S\colon \cQ \to \cY,\qquad \bq \coloneq (\bu,\bc) \mapsto S(\bq)
	= \by\coloneq (\rho,\phi,\bx)
\end{equation*}
with control and state spaces
\begin{align*}
	\cQ
	&\coloneqq \cU\times\cC :=
	L^\infty(0,T;\R^2)^M\times H^1(0,T)^M, \\
	\cY 
	&
	\coloneqq 
	W^{2,1}_p(Q_T)\times
	\paren[auto](){L^\infty(0,T;W_{\textup{DN}}^{2,p}(\Omega))\cap
	W^{1,p}(0,T;W^{1,p}(\Omega))}\times W^{1,s}(0,T;\R^2)^M
\end{align*}
for $s=\paren[auto](){\frac12-\frac1p}^{-1}$.
Furthermore we define the set of admissible controls
\begin{equation*}
	\Qad := \{(\bu,\bc)\in \cQ \colon \abs{u_i(t)} \le 1,\ 0\le c_i(t)\le 1\ 
	\text{f.a.a.\ }t\in (0,T)\ \text{and all}\ i=1,\ldots,M\}.
\end{equation*}

The optimal control problem we study in this article reads
\begin{subequations}\label{eq:optimal_control}
\begin{alignat}{3}
	&\text{Minimize}\quad & \cJ(\by; \bq) &\coloneqq \int_{\widetilde Q_T} e^{\nu\,t}\,\rho(t,x)\d x\d t
	-\mu\,\sum_{i=1}^M\int_0^T\ln(\xi(x_i(t)))\d t
	\nonumber\\
	&&&\qquad 
	+ \frac{\alpha_1}{2\,T}\sum_{i=1}^M\norm{u_i}_{H^1(0,T;\R^2)}^2
	+ \frac{\alpha_2}{2\,T}\sum_{i=1}^M \norm{c_i}_{H^1(0,T)}^2
	 \label{eq:optimal_control1}\\
	&\text{subject to}&
	\by&:=(\rho,\phi,\bx) = S(\bq),\\
	&& \bq &:= (\bu,\bc) \in \Qad.
\end{alignat}
\end{subequations}
The objective functional $\cJ$ aims at a fast 
evacuation of the crowd. By the factor $e^{\nu\,t}$ higher densities at a later time are penalized more. We observe the density in a subregion $\widetilde Q_T = I\times\widetilde\Omega$ where $\widetilde\Omega\subset \Omega$ is a subregion which the pedestrians must leave.
We use the temporal $H^1$-norm of the agent movement
directions and the intensities as a regularization to guarantee the smoothness required in \cref{thm:ex_state}. The regularization parameters $\alpha_1,\alpha_2>0$ are arbitrary but positive.

The fourth term in the objective is a barrier used to avoid that the
agents walk through walls. The barrier function $\xi\in H_D^1(\Omega)$
is the weak solution of the singularly perturbed problem
\begin{subequations}
\label{eq:barrier_function}
\begin{alignat}{3}
	-\delta_4\Delta \xi + \xi &= 1 &\qquad& \text{in}\ \Omega,\\
	\xi&= 0 && \text{on}\ \partial\Omega.
\end{alignat}
\end{subequations}
The barrier function $-\ln(\xi(x_i(t)))$ tends to infinity if $\dist(x_i(t),
\partial\Omega)\to 0$ for some $t\in [0,T]$. For $x_i(t)\in \interior{\Omega}$ 
there holds $\lim_{\mu\to 0} \mu\,\ln(\xi(x_i(t))) = 0$.
Here we choose $\mu >0$ to be fixed but small.

The control constraint $(\bu,\bc)\in \Qad$ guarantees that 
the agents do not move faster than the density in their current position
allows and that the intensity is bounded by reasonable values.

We have the following well-posedness result and necessary optimality condition.
\begin{theorem}
	\label{theorem:fonc}
	There exists at least one global solution $(\by,\bq)\in
	\cY\times \Qad$ of \eqref{eq:optimal_control}.

	Furthermore, each local minimizer $(\by,\bq)\in \cY\times \Qad$, $\by = (\rho,\phi,\bx)$, $\bq = (\bu,\bc)$,
	of \eqref{eq:optimal_control} fulfills for all directions in
	the tangential cone at $(\bu,\bc)$, namely $\delta \bq = (\delta \bu,\delta\bc) \in \cT_{\Qad}(\bu,\bc)$,
	\begin{align*}
		&\int_{\widetilde Q_T} e^{\nu\,t}\,\delta \rho(t,x)\d x\d t
		+ \frac{\alpha_1}{T} \,
		\inner{\bu}{\delta\bu}_{H^1(0,T;\R^2)^M}
		+ \frac{\alpha_2}{T} \,
		\inner{\bc}{\delta\bc}_{H^1(0,T)^M} \\
		&\qquad		
		-\mu\sum_{i=1}^M\int_0^T\frac{\nabla\xi(x_i(t))^\top\,\delta
		x_i(t)}{\xi(x_i(t))}\,\d t\ge 0,
	\end{align*}
	with $\by = S(\bu,\bc)$ and $\delta \by = (\delta \rho,\delta
	\phi,\delta \bx) = S'(\bu,\bc) \, (\delta\bu,\delta\bc)$ 
	characterized by the system
	\begin{subequations}
	\label{eq:linearized_system_strong}
	\begin{equation}\label{eq:linearized_system_strong_a}
		\partial_t  \delta\rho - \varepsilon \laplace \delta\rho - 
		\nabla\cdot\paren[auto](){\delta \rho \,\beta(\rho,\phi,\bx,\bc) + 
		\rho\paren[auto](){
		\frac{\partial\beta(\rho,\phi,\bx,\bc)}{\partial(\rho,\phi,\bx,\bc)}(\delta\rho,\delta\phi,\delta\bx,\delta \bc)
		}
		}
		= 
		0
		,
	\end{equation}
	\vspace{-4mm}
	\begin{alignat}{1}
		-\delta_1\, \laplace \delta\phi + 2\nabla\phi\cdot\nabla \delta\phi + \frac{2 f(\rho)\,f'(\rho)}{(f^2(\rho)+\delta_2)^2} \delta \rho
		&= 
		0 
		,
		\\
		\dot{\delta x}_i - v_0\,f'(\rho(\cdot,x_i))
		\left(		
		\nabla\rho(\cdot,x_i)^\transp \delta x_i
		+ \delta\rho(\cdot,x_i)
		\right)u_i 
		&=
		v_0\,f(\rho(\cdot,x_i)) \, \delta u_i,
	\end{alignat}
	\end{subequations}
	for $i=1,\ldots,M$, together with the boundary conditions \eqref{eq:forward_bc} and homogeneous initial conditions
	\begin{equation}
	\label{eq:linearized_system_boundary_conditions}
		\delta\rho(0,\cdot) 
		= 
		0
		\quad
		\text{and}
		\quad
		\delta x_i(0) = 0,
		\quad
		i=1,\ldots,M.
	\end{equation}
\end{theorem}
The proof of the theorem above is very close to those of Theorem 3.8 and Theorem 4.4 \cite{2011.03580}. The main difference is that the model in \cite{2011.03580} only allows to control the velocity $\bu$ of the agents, yet not their strength $\bc$. As for the existence proof, this does not impose any additional difficulty due to the uniform $L^\infty$-boundedness of $\bc$ as a consequence of the embedding $H^1 \hookrightarrow L^\infty$ in one spatial dimension. For the differentiability result, one has to add the derivatives with respect to $\bc$, yielding an additional term in \eqref{eq:linearized_system_strong_a} that, however, can be estimated similarly to the remaining terms.
\section{Discretization of the state equation}\label{sec:num}

In this section we introduce the numerical scheme used to compute approximate solutions of the forward
system \eqref{eq:forward_system}--\eqref{eq:forward_initial_conditions}.  To
this end, we introduce a semi-implicit time-stepping scheme and use a finite
volume discretization for the density function $\rho$ and continuous Lagrange
finite elements for the potential function.

\subsection{Space discretization}\label{sec:space_discretization}

For the spatial discretization of the system
\eqref{eq:forward_system}--\eqref{eq:forward_initial_conditions}
we define a family of geometrically conforming triangular meshes $\{\cT_h\}_{h>0}$.
For each $T\in \cT_h$ we denote by $h_T=\diam(T)$
the element diameter
and by $\rho_T$ the diameter of the largest inscribed ball in $T$. The
mesh parameter is then $h=\max_{T\in \cT_h} h_T$. The mesh family is
assumed to be shape regular meaning that there is a constant
$\kappa>0$ such that $h_T/\rho_T \le \kappa$ for all $T\in \cT_h$ and
all $h>0$. By $\cFinner$ we denote the set of interior element edges,
by $\cFboundary$ the boundary edges and write $\cF_h:= \cFinner\cup\cFboundary$.
Furthermore, to each edge $F\in\cF_h$ we associate the normal vector $n_F$
which is pointing outward in case of a boundary edge and has a fixed
orientation in case of an interior edge.

We propose a finite volume scheme for the transport equation. As we
use the finite element package \fenics\ for our implementation, we use a
notation which is rather usual for discontinuous Galerkin discretizations, see \cite{DiPietroErn2012} for an overview.
The finite-dimensional function spaces are defined by
\begin{align*}
 V_h &= \{ v \in L^\infty(\Omega)\colon
       	v|_T \in \mathcal{P}^0(T)\ \text{for all}\ T \in\mathcal{T}_h,\},\\
W_h &= \{ v \in C(\overline\Omega)\colon  \left. v \right|_T \in
\mathcal{P}^1(T) \ \mbox{for all}\ T \in\mathcal{T}_h\},\quad W_{h,\textup{D}} := W_h\cap H_D^1(\Omega),
\end{align*}
where $\mathcal{P}^k(T)$ denotes the space of polynomials on $T$ of degree not larger
than $k\in \N_0$. 
For a function $v\colon \Omega\to\R$, we define interface averages and jumps in the following way
\begin{equation*}
 \avg{v}_F := \frac12 (\left. v\right|_{T_1} + \left. v\right|_{T_2}),\qquad
 \jump{v}_F :=\left. v\right|_{T_1} - \left. v\right|_{T_2},\qquad \forall
 F\in \cFinner,
\end{equation*}
where $T_1, T_2\in\cT_h$ are chosen in such a way that $n_F = n_{\partial
T_1}{_{|F}} = - n_{\partial T_2}{_{|F}}$.

In order to discretize the system \eqref{eq:forward_system} we use
discontinuous approximations for $\rho$ and continuous ones for $\phi$.
The unknowns in our semi-discrete scheme are
\begin{equation*}
	\rho_h(t)\in V_h,\quad \phi_h(t)\in W_{h,\textup{D}},\quad x_1(t),\ldots,x_M(t)\in \R^2,\quad c_1(t),\ldots,c_M(t)\in \R
\end{equation*}
for all $t\in [0,T]$.
The approximate transport direction is then given by
\begin{equation*}
	\beta_h(\rho_h,\phi_h,\bx,\bc) \coloneq
	f(\rho_h)\,h(\nabla \phi_h + \phi_K(\bx,\bc))
\end{equation*}
with
\begin{equation*}
	\phi_K(\bx,\bc; t, x) := \sum_{j=1}^M c_j(t)\,K(x-x_i(t)).
\end{equation*}

The semi-discretization of \eqref{eq:forward_system1} then reads
\\[.5em]
Find\ $\rho_h\colon [0,T] \to V_h$ with $\rho_h(0) =
\text{proj}_{V_h}(\rho_0)$ and
\begin{equation}
\label{eq:semidiscrete_rho_equation}
	(\partial_t \rho_h(t), v_h)_\Omega + a(\rho_h(t),v_h) + b(\beta_h)(\rho_h(t),v_h) =
	0\qquad \forall v_h\in V_h, t\in (0,T).
\end{equation}
Here, $\text{proj}_{V_h}\colon L^2(\Omega)\to V_h$ is some projection operator,
$(\cdot,\cdot)_\Omega$ is the standard $L^2(\Omega)$-inner product and the bilinear forms are defined by
\begin{subequations}
\label{eq:bilinear_forms_rho}
\begin{align}
	a(\rho_h,v_h) &=  
	\varepsilon\sum_{F\in \cFinner} 
	\int_F \tau_F\,\jump{\rho_h}\jump{v_h} \d s	
	+ \sum_{F\in \cFboundary}
	\chi_{\partial\Omega_D}\gamma \int_F \rho_h\,v_h\d s
	\\	
	b(\beta_h)(\rho_h,v_h) &= - \sum_{F\in \cFinner} \int_F (\rho_h\,\beta_h)_F^*\,\jump{v_h}\,\d s.
\end{align}
\end{subequations}
The parameter $\tau_F$ is defined by
\begin{equation*}
	\tau_F \coloneq \frac{1}{\abs{x_{T_1}-x_{T_2}}},
\end{equation*}
where $x_{T}$ is the intersection of the orthogonal edge bisectors of $T\in
\cT_h$.
The term $\tau_F\,\jump{\rho_h}\jump{v_h}$ with $v_h=\chi_T$ for some $T\in
\cT_h$ approximates the diffusive flux
$\nabla\rho_h\cdot n_{\partial T}$ over the edge $F\subset T$. 
The bilinear form $b$ establishes the convective flux
$\beta\,\rho\,\cdot n_{\partial_T}$. As numerical flux function $(\cdot)^*$, we
choose the Lax-Friedrichs flux, see \cite{RiderLowrie2002}, defined by
\begin{equation}\label{eq:numerical_flux}
	(\rho_h\,\beta)_F^* = \avg{\rho_h\,\beta}_F\cdot n_F - \frac\eta2\,\jump{\rho_h}_F.
\end{equation}
The stabilization parameter $\eta\in \R$ is specified later.
For the closely related chemotaxis model such an approach has been sucessfully
applied in \cite{LiShuYang2017,GuoLiYang2018}. Of course, also other flux
functions are possible, e.g., the central upwind flux,
\cite{EpshteynKurganov2009}.

The Eikonal equation \eqref{eq:forward_system2} is discretized in
space using standard linear Lagrange elements which yields
\begin{equation}\label{eq:semidiscrete_Eikonal}
	\delta_1\,\paren(){\nabla \phi_h(t),\nabla w_h}_\Omega
	+ \paren[auto](){\abs{\nabla \phi_h(t)}^2,w_h}_\Omega =
	\paren[auto](){\frac{1}{f(\rho_h(t))^2+\delta_2},w_h}_\Omega
	\qquad\forall w_h\in W_{h,\textup{D}}.
\end{equation}
In our numerical experiments we used the Newton solver integrated in
\fenics. The Jacobian is established by automatic differentiation.

The ordinary differential equations for the agent trajectories
\eqref{eq:forward_system3} depend on a point evaluation
$\rho(t,x_i(t))$ of a function which is discontinuous in the discrete setting. In particular, this term would not be differentiable with respect to $x_i(t)$. As a remedy,  
we use instead of a point evaluation, see also \cite{Borsche:2015:1}, the following regularization
\begin{equation*}
	\rho_h(t,x_i(t)) \approx \eta_{x_i(t)}* \rho_h(t),
	\quad\text{with}\quad
	\eta_{x_i(t)}
	:=\frac{\delta_{x_i(t)}}{\delta_{x_i(t)}*1},
\end{equation*}
where $*$ stands for the convolution integral $\delta_{x_0}*v = \int_\Omega \delta_{x_0}\,v\,\d x$ of the functions $v\in L^1(\Omega)$ and some kernel function $\delta_{x_0}\in C^\infty(\R^2)$. An obvious choice is the regularized Dirac delta function
\begin{equation*}
	\delta_{x_0}(x) := \frac{1}{2\,\pi\,\zeta}\,\e^{-\frac{\norm{x-x_0}^2}{2\,\zeta}}
\end{equation*}
with small locality parameter $\zeta>0$.
Note that for $\zeta\to 0$ there holds $\delta_{x_0}*v \to v(x_0)$ for
any $v\in C(\overline\Omega)$.
Furthermore, the regularized Dirac delta fulfills $\int_{\R^2}
\delta_{x_0}\d x=1$ for arbitrary $x_0\in \R^2$.
The discretized ordinary differential equation then reads
\begin{equation}\label{eq:semidiscrete_ODE}
	\dot{x}_i(t) = 
	v_0\,f\paren[auto](){\eta_{x_i(t)}*\rho_h(t)}\,u_i(t),\quad n=1,\ldots,N
\end{equation}
and initial conditions $x_i(0) = x_{i,0}$.

\subsection{Time discretization}
For the temproal discretization we cover the time interval $[0,T]$ by an equidistant grid 
$I_\tau\coloneqq \{t_n\}_{n=0}^N$ with grid points
$t_n:= n\,\tau$ and grid size $\tau := T/N$. The spatial and temporal
discretization parameters are summarized in $\sigma = (h,\tau)$.
Moreover, we define the space of grid
functions
\begin{equation*}
	H_\tau(V) = \{v\colon I_\tau\to V\},
\end{equation*}
with $V$ an arbitrary linear space. If $V$ is again a function space
containing functions $v\colon \Omega\to \R$ we write $v(t_n) = v(t_n,\cdot)$.
The functions $\rho_h$, $\phi_h$, $x_i$, $u_i$ and $c_i$ arising in the semidiscrete
equations \eqref{eq:semidiscrete_rho_equation},
\eqref{eq:semidiscrete_Eikonal} and \eqref{eq:semidiscrete_ODE}
are approximated by grid functions
\begin{align*}
	&\rho_{\sigma}\in H_\tau(V_h),\quad \phi_{\sigma}\in H_\tau(W_{h,\textup{D}}),\quad
	u_{i,\sigma}, x_{i,\sigma}\in H_\tau(\R^2),\quad c_{i,\sigma}\in H_\tau(\R)
\end{align*}
for $i=1,\ldots,M$.
For brevity we write for $n=0,\ldots,N$
\begin{equation*}
	\rho_h^n \coloneq \rho_\sigma(t_n,\cdot), \quad
	\phi_h^n \coloneq \phi_\sigma(t_{n},\cdot), \quad
	x_i^n \coloneq x_{i,\sigma}(t_n), \quad
	u_i^n \coloneq u_{i,\sigma}(t_n)
\end{equation*}
and for the transport vector we use
\begin{equation*}
	\beta_h^n = \beta_h(\rho_h^n,\phi_h^n,\bx^n,\bc^n).
\end{equation*}

We replace the temporal derivative by a difference
quotient and use a semi-implicit time-stepping scheme, more precisely, the
diffusion-related terms are evaluated implicitly and the convection-related terms
explicitly. This yields the fully-discrete system
\begin{subequations}
\label{eq:full_discrete_forward}
\begin{align}
\label{eq:discrete_transport_eq}
	(\rho_h^{n+1},v_h)_\Omega + \tau\,a(\rho_h^{n+1}, v_h)
	&=
	(\rho_h^n, v_h)_\Omega - \tau\,b(\beta_h^{n})(\rho_h^n,v_h), \\
	\label{eq:discrete_eikonal_eq}
	\delta_1 \paren[auto](){\nabla\phi_h^n,\nabla w_h}_\Omega
	+ \paren[auto](){\abs{\nabla\phi_h^n}^2, w_h}_\Omega
	&= \paren[auto](){\frac{1}{f(\rho_h^n)^2+\delta_2}, w_h}_\Omega,\\
	\label{eq:discrete_agent_eq}
	x_i^{n+1} - x_i^n &=
	\tau\,v_0\,f\paren[auto](){\eta_{x^{n+1}_i}*\rho_h^{n+1}}\,u_i^{n+1},
\end{align}
\end{subequations}
for all test functions $v_h\in V_h$, $w_h\in W_{h,\textup{D}}$ and indices $i=1,\ldots,M$, $n=0,\ldots,N-1$. Furthermore, the initial conditions are established by means of:
\begin{equation*}
	\rho_h^0 = \proj{V_h}{(\rho_0)},\qquad x_i^0 = x_{i,0},\ i=1,\ldots,M.	
\end{equation*}
Note that the system of equations \eqref{eq:full_discrete_forward}
completely decouples and we can compute each variable after the other, in the following order
\begin{equation}\label{eq:computation_order}
	\rho_h^0, \bx^0 \mapsto \phi_h^0 \mapsto \rho_h^1 \mapsto \bx^1 \mapsto \phi_h^1 \mapsto\ldots \mapsto
	\rho_h^{N-1}\mapsto \bx^{N-1} \mapsto \phi_h^{N-1} \mapsto \rho_h^{N} \mapsto \bx^{N}.
\end{equation}

\subsection{Quality of discrete solutions}

In this section we study some basic properties for the solutions of
\eqref{eq:full_discrete_forward}.
In particular, it is of interest whether the physical bounds observed for the
solution of the continuous system
\eqref{eq:forward_system}--\eqref{eq:forward_initial_conditions} are
transferred to the discrete setting.

The basis functions of the finite element space $V_h$ are denoted by $\{\chi_T\}_{T\in \cT_h}$ defined by $\chi_T|_{T'} \equiv \delta_{T,T'}$ for all $T,T'\in \cT_h$. Note that by a slight abuse of notation we use the elements of $\cT_h$ as indices here. Introducing the matrices $M=(m_{T,T'})_{T,T'\in \cT_h}$, $A=(a_{T,T'})_{T,T'\in\cT_h}$ and $B^n=(b_{T,T'}^n)_{T,T'\in\cT_h}$ with entries
\begin{subequations}
\label{eq:matrix_entries}
\begin{align}
\label{eq:entries_mass_matrix}
	m_{T,T'} &= \begin{cases} \displaystyle |T|, &\text{if}\ T=T',\\ \displaystyle 0, &\text{otherwise},\end{cases} \\
\label{eq:entries_stiffness_matrix}
	a_{T,T'} &= \begin{cases}
	\displaystyle \varepsilon\,\sum_{F\in\cF_T\cap \cFinner}\tau_F\,\abs{F} + \gamma\sum_{F\in\cF_T\cap \cFboundary}\abs{F},
	&\text{if}\ T=T', \\
	\displaystyle -\varepsilon\,\tau_F\,\abs{F},
	&\text{if}\ T\ne T'\ \text{and}\ F:=\partial T\cap \partial T'\ne\emptyset, \\
	\displaystyle 0,& \text{otherwise},
	\end{cases} \\
\label{eq:entries_advection_matrix}	
	b_{T,T'}^n &= \begin{cases}
	\displaystyle -\frac12\sum_{F\in
	\cF_T\cap\cFinner}\paren[auto](){\int_{F}\beta_h^n|_T\cdot n_{\partial T}\d s
		   - \eta\,\abs{F}}, &\text{if}\ T=T',\\
	\displaystyle -\frac12 \paren[auto](){\int_F\,\beta_h^n|_{T'}\cdot n_{\partial T}\d s
	+ \eta\,\abs{F}},
	&\text{if}\ T\ne T'\ \text{and}\ F=\partial T\cap \partial T'\ne\emptyset, \\
	\displaystyle 0, &\text{otherwise},
	\end{cases}
\end{align}
\end{subequations}
allows to rewrite the system of equations \eqref{eq:discrete_transport_eq} as
\begin{equation}
\label{eq:matrix_vector_notation}
	(M+\tau\,A)\,\vec\rho^{n+1} = (M-\tau\,B^n)\,\vec\rho^n.
\end{equation}
Here, $\vec\rho^{n}$, $n=0,\ldots,N$, are the vector representations of $\rho_h^n$ with
respect to the basis $\{\chi_T\}_{T\in \cT_h}$. Note that the matrix $B^n$ depends also on $\vec\rho^n$.

\begin{theorem}
	The numerical scheme \eqref{eq:discrete_transport_eq} is mass conserving in the following sense. Assuming that $\gamma=0$ holds, i.\,e., there are no-flux boundary conditions present at all boundary parts $\partial\Omega_{\text{D}}$ and $\partial\Omega_{\text{W}}$, the solution $\rho_\sigma$ fulfills
	\begin{equation*}
		\int_\Omega \rho_h^{n}\d x = \int_\Omega \textup{proj}_{V_h}(\rho_0) \d x
		\quad\forall n=0,1,\ldots,N.
	\end{equation*}
\end{theorem}
\begin{proof}
	The assertion is trivially fulfilled for $n=0$ as the initial
	condition is established by $\rho_h^0=\textup{proj}_{V_h}(\rho_0)$.
	In matrix-vector notation the assertion is equivalent to $\vec
	1^\top M\vec\rho^{n+1} = \vec 1^\top M\vec\rho^{n}$.	
	This follows from \eqref{eq:matrix_vector_notation} after using
	\begin{equation*}
		\vec 1^\top A\vec\rho^{n+1} =
		\sum_{T\in \cT_h} \sum_{T'\in \cT_h}a_{T,T'}\rho_{T'}^{n+1} = 0
	\end{equation*}
	and	
	\begin{equation*}
		\vec 1^\top B^n\vec \rho^n =
		\sum_{T\in \cT_h} \sum_{T'\in \cT_h} b_{T,T'}\rho_{T'}^n.
	\end{equation*}
	In this expression, the stabilization
	terms (the ones multiplied with $\eta$) cancel out. Furthermore, after
	sorting terms in the sum by the edges
	$F\in \cFinner$ and denoting by $T_{F,1},T_{F,2}$ the two triangles meeting
	in $F$, we obtain the terms
	\begin{align*}
		\vec 1^\top B^n\vec \rho^n = -\frac12\sum_{F\in\cFinner}\int_F\Big(
		&  \beta_h|_{T_{F,1}}\,n_{\partial {T_{F,1}}}\,\rho_{T_{F,1}}
		+ \beta_h|_{{T_{F,2}}}\,n_{\partial {T_{F,1}}}\,\rho_{{T_{F,2}}} \\
		+ &\beta_h|_{{T_{F,2}}}\,n_{\partial {T_{F,2}}}\,\rho_{{T_{F,2}}}
		+ \beta_h|_{T_{F,1}}\,n_{\partial {T_{F,2}}}\,\rho_{T_{F,1}}\Big)\d s = 0.
	\end{align*}
	The last step follows due to $n_{\partial T_{F,2}} = - n_{\partial
	T_{F,2}}$ which implies $\vec 1^\top B^n\vec \rho^n = 0$.
\end{proof}

\begin{theorem}
Choose $\eta=1$ in \eqref{eq:numerical_flux} and denote by $\kappa>0$
the maximal aspect ratio of the mesh family $\cT_h$, see
\cref{sec:space_discretization}.
Let $\tau$ be chosen to satisfy the CFL condition 
\begin{equation}\label{eq:cfl_condition}
	\tau \le \frac{\pi}{3\,\kappa^2}\,\min_{T\in \cT_h} h_T.
\end{equation}
If furthermore, there holds $\proj_{V_h}(\rho_0) \in [0,1]$ a.\,e.\ in $\Omega$ and
$\beta_h^n= (1-\rho_h^n)\,\Phi^n$ with $\abs{\Phi^n} \le
1$, $n=0,\ldots,N$, the solutions of \eqref{eq:full_discrete_forward} fulfill for all $n=0,\ldots,N$
\begin{equation*}
	\rho_h^n(x)\in[0,1]\quad \text{f.a.a.\ } x\in \Omega.
\end{equation*}
\end{theorem}
\begin{proof}
	The diagonal entries of $(M+\tau\,A)$ are all positive
	and the off-diagonal entries are negative. Moreover, one easily
	concludes the strict diagonal dominance, this is,	
	\begin{equation*}
		\sum_{\genfrac{}{}{0pt}{}{T'\in\cT_h}{T'\ne T}}\abs{m_{T,T'}+\tau\,a_{T,T'}}
		< m_{T,T}+\tau\,a_{T,T}.
	\end{equation*}
	This implies that $(M+\tau\,A)$ is an M-matrix and
	consequently, the inverse $(M+\tau\,A)^{-1}$ exists and fulfills $(M+\tau\,A)^{-1} \ge 0$.
	
	Let $n\in \N_0$ be fixed and assume that $\rho_h^n(x)\in [0,1]$ for
	almost all $x\in \Omega$. We show $\rho_h^{n+1}\ge 0$ 
	by confirming that the right-hand side of
	\eqref{eq:matrix_vector_notation} has non-negative entries only.
	Assuming that $\rho_T^n\ge 0$, $T\in \cT_h$,
	we show that the entries of $(M-\tau\,B^n)$ are non-negative as well.
	The entries on the diagnoal have the form
	\begin{equation*}
		\abs{T}
		   + \frac\tau2\sum_{F\in \cF_T\cap \cFinner}\int_{F}(\beta_h^n|_T\cdot n_{\partial T}-1)\d s \ge \abs{T} - \tau\,\abs{\partial T},
	\end{equation*}
	where the first step follows from the assumption \eqref{eq:def_beta}
	implying $\abs{\beta_h^n} \le 1$.
	To estimate further we take the geometric mesh quantities and
	the shape regularity $h_T \le \kappa\,\rho_T$ into
	account and arrive at
	\begin{equation*}
		\frac{\abs{T}}{\abs{\partial T}}
		\ge \frac{\pi\,\rho_T^2}{3\,h_T} \ge
		\frac{\pi\,h_T}{3\,\kappa^2}
		\ge \tau,
	\end{equation*}
	where the last step is the CFL condition
	\eqref{eq:cfl_condition}. The previous two estimates confirm
	$(m_{T,T}-\tau\,b^n_{T,T}) \ge 0$ for all $T\in \cT_h$.
	The remaining entries in the matrix, namely $(m_{T,T'} -
	\tau\,b_{T,T'}^n)$ with $F=\partial T\cap \partial T'\ne\emptyset$,
	have the form
	\begin{equation*}
		\frac\tau2\,\int_F(\beta_h^n|_{T'}\cdot n_{\partial T}+1)\d s \ge 0.
	\end{equation*}
	The non-negativity follows again from $\abs{\beta_h^n}\le
	1$. Combining the previous arguments provides the lower bound $\rho_h^{n+1}\ge 0$.
	
	To show the upper bound we rearrange the equation system \eqref{eq:matrix_vector_notation}
	in the form
	\begin{equation}
	\label{eq:matrix_vector_notation_upper_bound}
		(M+\tau\,A)\,(\vec 1-\vec\rho^{n+1}) = M\,(\vec 1-\vec\rho^n)
		+ \tau\,B^n\,\vec\rho^n + \tau\,A\,\vec 1.
	\end{equation}
	We may rewrite the transport term using
	$\rho_h^n\,\beta_h^n=(1-\rho_h^n)\,\widetilde\beta_h^n$
	with $\widetilde\beta_h^n = \rho_h^n\,\Phi$.
	In \eqref{eq:matrix_vector_notation_upper_bound} we reformulate the
	expression involving $B^n$ by means of
	\begin{align*}
		[B^n\,\vec\rho^n]_T
		&= -\sum_{F\in \cF_T\cap \cFinner}
		\int_F\paren[auto](){\avg{(1-\rho_T^n)\,\widetilde\beta_h^n}_F\cdot n_F
		+ \frac12\,\jump{1-\rho_h^n}}_F\jump{\chi_T}_F\d s \\
		&= -\frac12 \sum_{F\in \cF_T\cap \cFinner}\int_F
		(\widetilde\beta_h^n|_T \cdot n_{\partial T}+1)\d s\cdot
		(1-\rho_T^n) \\
		&\qquad- \frac12\sum_{F\in \cF_T\cap \cFinner}
		\int_F(\widetilde\beta_h^n|_{T_F}\cdot n_{\partial
		T}-1)\d s\cdot (1-\rho_{T_F}^n) \\
		&=: [\widetilde B^n(1-\vec\rho^n)]_T.
	\end{align*}
	with $T_F\in \cT_h\setminus\{T\}$, $T_F\cap T=F$.
	With the same arguments like above one can show that the
	entries of $M+\tau\,\widetilde B^n$ are non-negative and
	together with $1-\vec\rho^n\ge 0$, $A\,\vec 1\ge 0$ and the
	M-matrix property of $M+\tau\,A$ we
	arrive at the desired bound $1-\vec\rho^{n+1}\ge 0$.
\end{proof}

\subsection{The discrete optimal control problem}

Next, we study a discrete version of the optimal control problem
\eqref{eq:optimal_control}. The control and state variables are grid functions in time and thus, we introduce the discrete $H^1(0,T)$-inner product for functions $u_\tau, v_\tau\colon \{t_n\}_{n=0}^N\to V$, with $V$ some Hilbert space,
\begin{equation*}
	(u_\tau,v_\tau)_{H^1(0,T;V),\tau} := \tau \sum_{n=0}^{N} \paren[auto](){u^n, v^n}_{V\times V} +
	\tau^{-1} \sum_{n=0}^{N-1} \paren[auto](){u^{n+1} - u^n,v^{n+1} - v^n}_{V\times V}.
\end{equation*}
This induces the norm $\norm{u_\tau}_{H^1(0,T;V),h}^2 :=
(u_\tau,u_\tau)_{H^1(0,T;V),h}$.
For the discrete control space we obtain
\begin{align*}
	\cU_\sigma &:=
	\{
	\bu_\sigma = (u_{1,\sigma},\ldots,u_{M,\sigma})\colon
	\ u_{i,\sigma} \in H_\tau(\R^2)\ \text{for}\ i=1,\ldots,M\}
	\} \\	       
	\cC_\sigma &:= \{\bc_\sigma = (c_{1,\sigma},\ldots,c_{M,\sigma})\colon\
	c_{i,\sigma} \in H_\tau(\R)\ \text{for}\ i=1,\ldots,M\},\\
	\cQ_\sigma &:= \cU_\sigma\times\cC_\sigma,
\end{align*}
and the admissible set by
\begin{align*}
	\cU_{\sigma,\textup{ad}} &:= \{\bu_\sigma\in \cU_\sigma\colon \abs{u_i^n}\le 1,\ i=1,\ldots,M,\ n=0,\ldots,N\}, \\
	\cC_{\sigma,\textup{ad}} &:= \{\bc_\sigma\in \cC_\sigma\colon 0\le c_i^n \le 1, i=1,\ldots,M, n=0,\ldots,N\}, \\
	\cQ_{\sigma,\textup{ad}} &:= \cU_{\sigma,\textup{ad}}\times \cC_{\sigma,\textup{ad}}.
\end{align*}
The discrete state space is defined by
\begin{equation*}
	\cY_\sigma := H_\tau(V_h)\times H_\tau(W_{h,\textup{D}})\times H_\tau(\R^2)^M.
\end{equation*}
With these definitions, the discrete optimal control problem related to
\eqref{eq:optimal_control} reads as 
\begin{subequations}\label{eq:discrete_ocp}
\begin{alignat}{2}
	&\text{Minimize} &\quad
	\cJ_\sigma(\by_\sigma,\bq_\sigma) &= \tau\,\sum_{n=1}^N e^{\nu\, t_n}\,\int_\Omega\,\rho_h^n(x)\d x
	-\mu\,\tau\sum_{i=1}^M\sum_{n=1}^N \ln(\eta_{x_i^n}*\xi_h) \nonumber\\
	&&&\quad + \frac{\alpha_1}{2\,T} \sum_{i=1}^M \norm{u_{i,\sigma}}_{H^1(0,T),\tau}^2
	+ \frac{\alpha_2}{2\,T} \sum_{i=1}^M \norm{c_{i,\sigma}}_{H^1(0,T),\tau}^2	
	\label{eq:discrete_ocp_objective}
	\\
	&\text{subject to} &
	&\by_\sigma:=(\rho_\sigma,\phi_\sigma,\bx_\sigma) = S_\sigma(\bq_\sigma), \\
	&&& \bq_\sigma := (\bu_\sigma,\bc_\sigma)\in \cQ_{\sigma,\textup{ad}},
\end{alignat}
\end{subequations}
where is $\xi_h\in W_h$ the finite element approximation of
\eqref{eq:barrier_function} with first-order Lagrange elements.
Furthermore, $S_\sigma$ is the solution operator of
\eqref{eq:full_discrete_forward}. Note that, in order to maintain the
differentiability of the barrier term with respect to $x_i^n$, we use
a regularization of the point evaluation of the nonsmooth function $\xi_h$, 
compare also \eqref{eq:discrete_agent_eq}.

We may write the control problem \eqref{eq:discrete_ocp} in the more
compact reduced form 
\begin{equation}
\label{eq:reduced_ocp}
	j_\sigma(\bq_\sigma) := \cJ_\sigma(S_\sigma(\bq_\sigma),\bq_\sigma)\to
	\min!\quad\text{subject to}\quad
	\bq_\sigma\in \cQ_{\sigma,\textup{ad}}.
\end{equation}

To deduce a necessary optimality condition we apply the Lagrange formalism. 
The Lagrange function 
\begin{equation*}
	\cL_\sigma\colon \cY_\sigma\times \cQ_\sigma \times \cY_\sigma
	\to \R
\end{equation*}
coupling the discrete state equation
\eqref{eq:full_discrete_forward} reads
\begin{align*}
	 &\cL_\sigma(\rho_\sigma,\phi_\sigma,\bx_\sigma; \bu_\sigma, \bc_\sigma; \lambda_{\rho,\sigma}, \lambda_{\phi,\sigma},\lambda_{\bx,\sigma})	 
	  = \cJ_{\sigma}(\rho_\sigma,\phi_\sigma,\bx_\sigma;
	  \bu_\sigma, \bc_\sigma) \\
	  &\quad
	  -\int_\Omega(\rho_h^0-\proj_{V_h}(\rho_0))\,\lambda_{\rho,h}^0\d x
	  - \sum_{n=0}^{N-1} \paren[auto](){\int_\Omega(\rho_h^{n+1} -
	  \rho_h^{n})\,\lambda_{\rho,h}^{n+1}\d x
	  + \tau\,a(\rho_h^{n+1},\lambda_{\rho,h}^{n+1}) -
	  \tau\,b(\beta_h^n)(\rho_h^n,\lambda_{\rho,h}^{n+1})}
	  \\
	  &\quad- \sum_{n=0}^{N-1} \tau \paren[auto](){\delta_1\int_\Omega
	  \nabla \phi_h^{n}\cdot\nabla \lambda_{\phi,h}^n \d x
	  + \int_\Omega \abs{\nabla \phi_h^n}^2\,\lambda_{\phi,h}^n\d x - \int_\Omega \frac{1}{f(\rho_h^n)^2 + \delta_2}\,\lambda_{\phi,h}^n\d x} \\
	  &\quad- \sum_{i=1}^M\paren[auto](){
	  (x_i^0-x_{i,0})^\top\,\lambda_{\bx,i}^0 +
	  \sum_{n=0}^{N-1}\paren[auto](){ x_i^{n+1} - x_i^n -
	  \tau\,f(\eta_{x_i^{n+1}}*\rho_h^{n+1})\,u_i^{n+1}}^\top\,\lambda_{\bx,i}^{n+1}	  
	  }.	  
\end{align*}
To shorten the notation we write $\blambda_\sigma:=
(\lambda_{\rho,\sigma},\lambda_{\phi,\sigma},\lambda_{\bx,\sigma})$.
The adjoint equation system determining these variables
for a given control and state is
\begin{subequations}
\label{eq:discrete_adjoint_system}
\begin{alignat}{6}
	&\lambda_{\rho,h}^{n}\in V_h\colon&&
	\frac{\partial\cL_\sigma}{\partial\rho_h^n}(\by_\sigma, \bq_\sigma,\blambda_{\sigma})\delta\rho_h &= 0
	&&&\forall \delta\rho_h\in V_h,\ n=0,\ldots,N, \\
	&\lambda_{\phi,h}^{n}\in W_{h,\textup{D}}\colon&\quad&
	\frac{\partial\cL_\sigma}{\partial\phi_h^n}(\by_\sigma, \bq_\sigma,\blambda_{\sigma})\delta\phi_h &= 0
	&&&\forall \delta\phi_h\in W_h,\ n=0,\ldots,N-1, \\
	&\lambda_{\bx_i}^{n}\in \R^2 \colon&&
	\frac{\partial\cL_\sigma}{\partial x_i^n}(\by_\sigma,\bq_\sigma,\blambda_{\sigma})\delta x_i &= 0
	&\quad&&\forall \delta x_i\in \R^2,\ n=0,\ldots,N
\end{alignat}
\end{subequations}
for $i=1,\ldots,M$.
Note that this can be interpreted as a coupled system involving a parabolic
PDE and an ODE that run backward in time. We use the automatic differentiation
feature in \fenics\ in our implementation.
As the forward system completely decouples in each time step, so does the
adjoint system and we can compute step by step:
\begin{equation*}
	\lambda_{\bx,h}^{N} \mapsto \lambda_{\rho,h}^{N} \mapsto (\lambda_{\phi,h}^{N})\mapsto \ldots \mapsto
	\lambda_{\bx,h}^{0} \mapsto \lambda_{\rho,h}^{0} \mapsto \lambda_{\phi,h}^0.
\end{equation*}
With the adjoint states at hand we can assemble the derivatives of the
reduced objective \eqref{eq:reduced_ocp} and end up with the following
optimality condition for \eqref{eq:discrete_ocp}:
\begin{theorem}[Necessary optimality condition]
Let $(\by_\sigma,\bq_\sigma)\in \cY_\sigma\times \cQ_{\sigma,\textup{ad}}$ be a local solution of
\eqref{eq:discrete_ocp}. Then, there exists $\blambda_\sigma\in
\cY_\sigma$ fulfilling
\eqref{eq:discrete_adjoint_system} and
\begin{align}
\label{eq:optimality_condition_reduced_objective}
	&\alpha_1\,
	(u_{\sigma,i},v_{\sigma,i}-u_{\sigma,i})_{H^1(0,T;\R^2),\tau}
	+
	\alpha_2\,(c_{\sigma,i}, d_{\sigma,i}-c_{\sigma,i})_{H^1(0,T),\tau} \nonumber\\
	&\qquad
	+
	\tau\,\sum_{n=1}^{N}
	\paren[auto](){f(\eta_{x_i^n}*\rho_h^{n})\,\lambda_{\bx_i}^n,v_i^n-
	u_i{^n}}_{\R^2}
	+ \tau\,\sum_{n=0}^{N-1} \frac{\partial
	\paren[auto](){b(\beta_h^n)(\rho_h^n,\lambda_{\rho,h}^{n+1})}}{\partial
	c_i^n}(d_i^n-c_i^n)
	\ge 0 
\end{align}
for all test functions $\br_\sigma:=(\bv_\sigma,\bd_\sigma)\in
\cQ_{\sigma,\textup{ad}}$ and all $i=1,\ldots,M$.
\end{theorem}
\begin{proof}
It is well-known that the variational inequality
$j_\sigma'(\bq_\sigma)(\br_\sigma-\bq_\sigma)\ge 0$ for $\br_\sigma\in \cQ_{\sigma,\textup{ad}}$ is
necessary for $\bq_\sigma$ being a local minimizer of \eqref{eq:reduced_ocp}. Taking into account the the equivalence
\begin{equation*}
	j_\sigma'(\bq_\sigma)\delta\bq_\sigma =
	\frac{\partial\cL}{\partial\bq_\sigma}(\by_\sigma,\bq_\sigma,\blambda_\sigma)\delta\bq_\sigma
	\quad\text{if} 
	\blambda_\sigma\ \text{solves}\ \text{\eqref{eq:discrete_adjoint_system}}
\end{equation*}
yields the variational inequality \eqref{eq:optimality_condition_reduced_objective}.
\end{proof}

Our solution algorithm is based on a projected gradient algorithm and it remains to establish a representation of the gradient of $j_\sigma$.

The derivative of the objective \eqref{eq:discrete_ocp_objective} towards some direction $\delta \bq_\sigma= (\delta \bu_\sigma,\delta\bc_\sigma)\in \cQ$ reads
\begin{align*}
	j_\sigma'(\bq_\sigma)\delta \bq_\sigma = \sum_{i=1}^M\Big(&\alpha_1 (u_{i,\sigma}, \delta u_{i,\sigma})_{H^1(0,T),\tau} + \alpha_2 (c_{i,\sigma}, \delta c_{i,\sigma})_{H^1(0,T),\tau} \\
	&+ \tau \sum_{n=1}^{N} f(\eta_{x_i^n}*\rho_h^{n})\,\delta {u_i^n}^\top \lambda_{\bx_i}^n
	+ \tau\,\sum_{n=0}^{N-1} \frac{\partial
	\paren[auto](){b(\beta_h^n)(\rho_h^n,\lambda_{\rho,h}^{n+1})}}{\partial
	c_i^n} \delta c_i^n\Big).
\end{align*}
To obtain a representation of the $H^1(0,T),\tau$-gradient of $j_\sigma$ with respect to $\bu_\sigma$ and $\bc_\sigma$, we introduce the grid functions $z_{i,\sigma}: \{t_n\}_{n=0}^{N}\to \mathbb R^2$ and $d_{i,\sigma}:\{t_k\}_{n=0}^N\to\R$, $i=1,\ldots,M$ solving
\begin{subequations}
\label{eq:system_h1_representation}
\begin{align}
	&\left(
	\frac1{\tau^2}
	\begin{pmatrix}
	1 & -1 & & &\\
	-1 & 2 & -1 & &\\
	& -1 & 2 & -1 & \\	
	& \ddots & \ddots & \ddots & \\
	& & -1 & 2 & -1 \\
	& & & -1 & 1	
	\end{pmatrix}
	+
	I_{N+1\times N+1}
	\right)
	\begin{pmatrix}
		z_{i}^0 \\
		z_{i}^1 \\
		z_i^2 \\
		\vdots \\
		z_{i}^{N-1} \\
		z_{i}^{N}
	\end{pmatrix}
	=
	\begin{pmatrix}
	0 \\	
	-f(\rho_h^{1}(\bx_i^{1}))\lambda_{x_i}^0 \\
	-f(\rho_h^{2}(\bx_i^{2}))\lambda_{x_i}^1 \\
	\vdots \\
	-f(\rho_h^{N-1}(\bx_i^{N-1}))\lambda_{x_i}^{N-2}\\
	-f(\rho_h^{N}(\bx_i^{N}))\lambda_{x_i}^{N-1}
	\end{pmatrix}	
\end{align}
and
\begin{align}
	&\left(
	\frac1{\tau^2}
	\begin{pmatrix}
	1 & -1 & & &\\
	-1 & 2 & -1 & &\\
	& -1 & 2 & -1 & \\	
	& \ddots & \ddots & \ddots & \\
	& & -1 & 2 & -1 \\
	& & & -1 & 1	
	\end{pmatrix}
	+
	I_{N+1\times N+1}
	\right)
	\begin{pmatrix}
		d_{i}^0 \\
		d_{i}^1 \\
		d_i^2 \\
		\vdots \\
		d_{i}^{N-1} \\
		d_{i}^{N}
	\end{pmatrix}
	=
	\begin{pmatrix}
	-\partial_{c_i^0} b(\beta_h^0)(\rho_h^0,\lambda_{\rho,h}^1) \\
	-\partial_{c_i^1} b(\beta_h^1)(\rho_h^1,\lambda_{\rho,h}^2) \\
	-\partial_{c_i^2} b(\beta_h^2)(\rho_h^2,\lambda_{\rho,h}^3) \\
	\vdots \\
	-\partial_{c_i^{N-2}} b(\beta_h^{N-2})(\rho_h^{N-2},\lambda_{\rho,h}^{N-1}) \\
	-\partial_{c_i^{N-1}} b(\beta_h^{N-1})(\rho_h^{N-1},\lambda_{\rho,h}^{N}) \\
	0 \\
	\end{pmatrix}	
\end{align}
\end{subequations}

for $i=1,\ldots,M$. By a simple calculation we then confirm
\begin{align*}
	(z_{i,\sigma},\delta u_{i,\sigma})_{H^1(0,T; \R^2),\tau} &=
	\sum_{n=1}^{N} f(\eta_{x_i^{n}} * \rho_h^{n})\,\delta {u_i^n}^\top
	\lambda_{\bx_i}^n, \\
	(d_{i,\sigma},\delta c_{i,\sigma})_{H^1(0,T),\tau} &=
	\sum_{n=0}^{N-1} \frac{\partial\paren[auto](){b(\beta_h^n)(\rho_h^n,\lambda_{\rho,h}^{n+1})}}{\partial
	c_i^n} \delta c_i^n
\end{align*}
We write $\bz_\sigma = (z_{1,\sigma},\ldots,z_{M,\sigma})\in \cU$ and $\bd_\sigma = (d_{1,\sigma},\ldots,d_{M,\sigma})\in \cC$ and get the following representation of the gradient of $j_\sigma$:
\begin{subequations}
\label{eq:discrete_gradient}
\begin{align}
	\nabla_{\bu_\sigma} j_\sigma(\bq_\sigma) = \alpha_1\,\bu_{\sigma} + \bz_{\sigma}, \\
	\nabla_{\bc_\sigma} j_\sigma(\bq_\sigma) = \alpha_2\,\bc_{\sigma} + \bd_{\sigma}.
\end{align}
\end{subequations}
This allows an implementation of a projected gradient method
which we discuss in the following section.

\subsection{Optimization algorithms for the discretized problem}\label{sec:num_opt}

For a solution of the discretized optimal control problem
\ref{eq:discrete_ocp} we propose a \emph{projected gradient
algorithm}. In this procedure, for a given initial control $\bq^{(0)} = (\bu^{(0)},\bc^{(0)}) \in \cQ$, the new iterates are sucessively computed by means of
\begin{subequations}
\label{eq:projected_gradient_iteration}
\begin{align}
	\bu_\sigma^{(k+1)} = \bPi_{\textup{ad}}^u\paren[auto](){\bu_\sigma^{(k)} -
	s_{(k)}\,\nabla_{\bu_\sigma} j_\sigma(\bq_\sigma^{(k)})},\\
	\bc_\sigma^{(k+1)} = \bPi_{\textup{ad}}^c\paren[auto](){\bc_\sigma^{(k)} -
	s_{(k)}\,\nabla_{\bc_\sigma} j_\sigma(\bq_\sigma^{(k)})},	
\end{align}
\end{subequations}
with $\bPi_{\textup{ad}}^u\colon \cU_\sigma\to \cU_{\sigma,\textup{ad}}$ and $\bPi_{\textup{ad}}^c:\cC_\sigma\to\cC_{\sigma,\textup{ad}}$ the $H^1(0,T),\tau$ projections onto the admissible sets $\cU_{\sigma,\textup{ad}}$ and $\cC_{\sigma,\textup{ad}}$, respectively, this is,
\begin{subequations}
\label{eq:projection_problem}
\begin{align}
	\bPi_{\textup{ad}}^u(\bu_\sigma):= \argmin_{\bv_\sigma \in \cU_{\sigma,\textup{ad}}}
	\frac12\norm{\bu_{\sigma} - \bv_{\sigma}}_{H^1(0,T;\R^2),\tau}^2,\label{eq:projection_problem_u}\\
	\bPi_{\textup{ad}}^c(\bc_\sigma):= \argmin_{\bd_\sigma \in \cC_{\sigma,\textup{ad}}}
	\frac12\norm{\bc_{\sigma} - \bd_{\sigma}}_{H^1(0,T),\tau}^2.
\end{align}
\end{subequations}

A formula for the gradient of $j_\sigma$ has been derived in the previous
section already, see \eqref{eq:discrete_gradient}.
The step length parameter $s^{(k)}>0$ is obtained by an Amijo line search and
must fulfill the sufficient decrease condition
\begin{equation}\label{eq:armijo_interpolation}
	j_\sigma(\bq_\sigma^{(k)}-s_{(k)}\,\nabla j_\sigma(\bq_\sigma^{(k)})) \le
	j_\sigma(\bq_\sigma^{(k)}) - \frac{d}{s^{(k)}} \norm{\bq_\sigma^{(k)} - \paren[auto](){\bq_\sigma^{(k)}-s^{(k)}\,\nabla j_\sigma(\bq_\sigma^{(k)})}}_{\cQ_\sigma}^2
\end{equation}
with a decrease parameter $d\in (0,1)$ which is usually small (e.g. $10^{-4}$).
A reasonable stopping criterion for the projected gradient algorithm is
\begin{equation*}
	\norm{\bq_\sigma^{(k)} - \bPi_{\textup{ad}}\paren[auto](){\bq_\sigma^{(k)}
	- \nabla j_\sigma(\bq_\sigma^{(k)})}}_{\cQ_\sigma} \le 10^{-3}.
\end{equation*}

It remains to discuss the realization of the projection operators and we propose a primal dual active set strategy that may also be considered as semismooth Newton method. Note that the operators $\Pi_{\textup{ad}}$ are semismooth, see \cite{2112.12018}.
The evaluation of the projection operator $\bPi_{\textup{ad}}^u\colon \cU_\sigma\to
\cU_{\sigma,\textup{ad}}$ requires to solve the optimization problem \eqref{eq:projection_problem_u}.
The unknowns (assuming $M=1$ and omitting the agent's index $i$ for a while) are the coefficients of the functions $\cU_\sigma\ni \bPi_{\textup{ad}}(u_\sigma) = w_\sigma \simeq \vec w\in \R^{(N+1)\times 2}$ for some given $\cU_\sigma\ni u_\sigma \simeq \vec u\in \R^{(N+1)\times 2}$. We switch
to a matrix-vector notation and define
\begin{equation*}
	w^n := \begin{pmatrix}w^n_1\\  w_2^n\end{pmatrix} := w_\sigma(t_n) ,\qquad 
	\vec w_j := (w^0_j, \ldots, w^N_j)^\top,\ j=1,2,
\end{equation*}
as well as the matrix $A\in \R^{(N+1)\times(N+1)}$ on the left-hand side of the linear system
\eqref{eq:system_h1_representation} inducing the discrete
$H^1(0,T),\tau$-norm.
The Lagrangian for \eqref{eq:projection_problem} reads
\begin{equation*}
	L(\vec w_1, \vec w_2, \vec\lambda) 
	= \frac12\sum_{i=1}^2\paren[auto](){\vec w_i-\vec u_i}^\top
	A\paren[auto](){\vec w_i-\vec u_i}
	- \frac12\,\lambda^\top \paren[auto](){\abs{\vec w}_*^2 - \vec 1},
\end{equation*}
with $\abs{\vec w}_*^2 = (\abs{w^0}^2,\ldots,\abs{w^N}^2)^\top$. The
Karush-Kuhn-Tucker system for \eqref{eq:projection_problem} then reads
\begin{align*}
	A\,(\vec w_i-\vec u_i) - \vec\lambda\cdot \vec w_i = 0\qquad i=1,2,&\\
	\frac12\,\paren[auto](){\abs{\vec w}_*^2-\vec 1} \le 0,\qquad \vec\lambda
	\ge 0,\qquad \frac12\,\vec\lambda\cdot \paren[auto](){\abs{\vec
	w}_*^2-1} &= 0,
\end{align*}
where $\cdot$ is the component-wise multiplication of two vectors.
We reformulate the complementarity condition by means of a nonsmooth equation
and arrive at the following equivalent form of the KKT system
\begin{equation}
\label{eq:nonlinear_system_projection}
	F(\vec w_1,\vec w_2,\vec \lambda)
	:=
	\begin{bmatrix}
	A\,(\vec w_1-\vec u_1) - \vec\lambda\cdot \vec w_1\\
	A\,(\vec w_2-\vec u_2) - \vec\lambda\cdot \vec w_2\\
	\vec\lambda-\max\{0, - \frac12\,(\abs{\vec w}_*^2-\vec 1) + \vec\lambda\}
	\end{bmatrix}
	=
	0.
\end{equation}
This nonlinear system can be solved iteratively by a semismooth Newton
method. Given is an initial pair $(\vec u^{(0)},
\vec\lambda^{(0)})$. Successively, one computes the active and inactive set
\begin{align*}
	\cA^{(k)} &:=
	\{n\in \{0,\ldots,N\}\colon - \frac12\,(\abs{w^n}_2^2-1) +
	\lambda^n > 0\}, \\
	\cI^{(k)} &:= \{0,\ldots,N\}\setminus \cA^{(k)},
\end{align*}
solves the Newton system
\begin{equation*}
	\begin{bmatrix}
		A-D_{\vec\lambda^{(k)}} & 0 & -D_{\vec w_1^{(k)}} \\
		0 & A-D_{\vec\lambda^{(k)}} & -D_{\vec w_2^{(k)}} \\
		D_{\cA^{(k)}}\,D_{\vec w_1^{(k)}} &
		D_{\cA^{(k)}}\,D_{\vec w_2^{(k)}} &
		D_{\cI^{(k)}}
	\end{bmatrix}
	\begin{bmatrix}
		\vec{\delta w}_1 \\ \vec{\delta w}_2 \\ \vec{\delta \lambda}
	\end{bmatrix}
	=
	-
	\begin{bmatrix}
		A\,(\vec w_1^{(k)}-\vec u_1) - \vec\lambda^{(k)}\cdot \vec
		w_1^{(k)}\\ 
		A\,(\vec w_2^{(k)}-\vec u_2) - \vec\lambda^{(k)}\cdot \vec
		w_2^{(k)} \\
		\vec\lambda^{(k)}-\max\{0, - \frac12\,(\abs{\vec w^{(k)}}_*^2-\vec 1) + \vec\lambda^{(k)}\},
	\end{bmatrix}
\end{equation*}
with the diagonal matrices $D_{\vec v} = \text{diag}(\vec v)$ for $\vec v\in
\R^{N+1}$ and $D_{\cM} = \text{diag}(\chi_{\cM})$ for
$\cM\subset\{0,\ldots,N\}$, and performs the Newton update
\begin{equation*}
	\vec w^{(k+1)} = \vec w^{(k)} + \vec{\delta w},\qquad
	\vec \lambda^{(k+1)} = \vec \lambda^{(k)} + \vec{\delta \lambda}.
\end{equation*}
This procedure is repeated for $k=0,1,\ldots$ until some termination
criterion, e.\,g., $\norm{F(\vec w_1,\vec w_2,\vec \lambda)} < \text{tol}$, is
fulfilled.

\section{Numerical experiments}\label{sec:num_ex}

This section is devoted to numerical experiments.
To establish the discretized system \eqref{eq:discrete_ocp} the finite
element library \fenics\ was used, complemented by a \python\
implementation of the projected gradient method from
\cref{eq:projected_gradient_iteration} and the  Armijo step size rule from
\eqref{eq:armijo_interpolation}. The computational meshes were created by the mesh
generator \texttt{mshr} integrated in \fenics.

\subsection{Example 1}
\label{sec:example_1}

\begin{figure}
\subfloat[Solution at time $t=0$\label{fig:example_1_initial}]{
\centering
\fbox{\includegraphics[width=.49\textwidth]{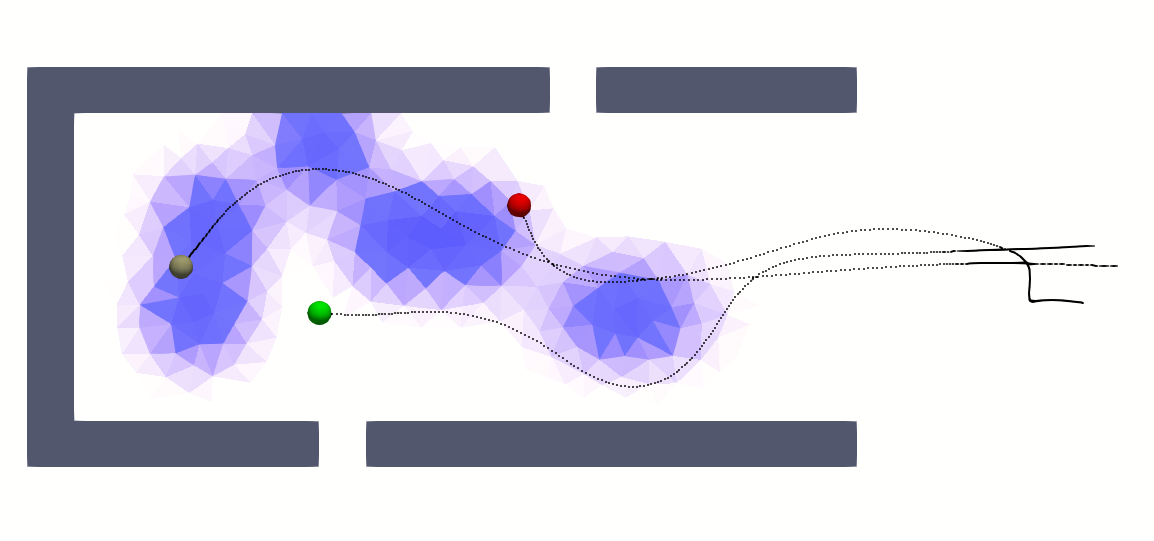}}
}
\subfloat[Solution at time $t=1.2$]{
\centering
\fbox{\includegraphics[width=.49\textwidth]{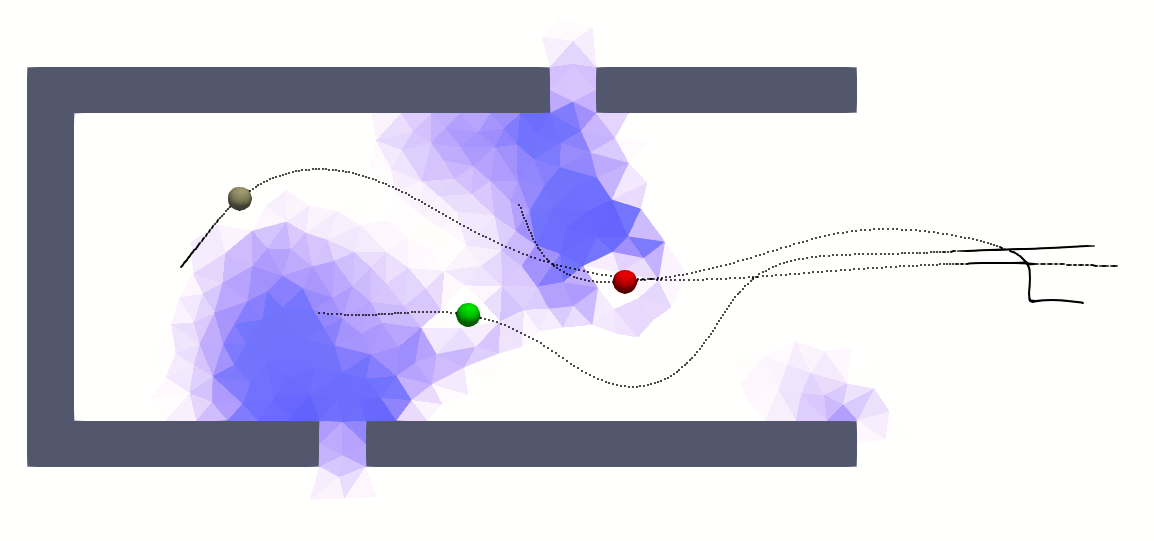}}
}

\subfloat[Solution at time $t=2.4$]{
\centering
\fbox{\includegraphics[width=.49\textwidth]{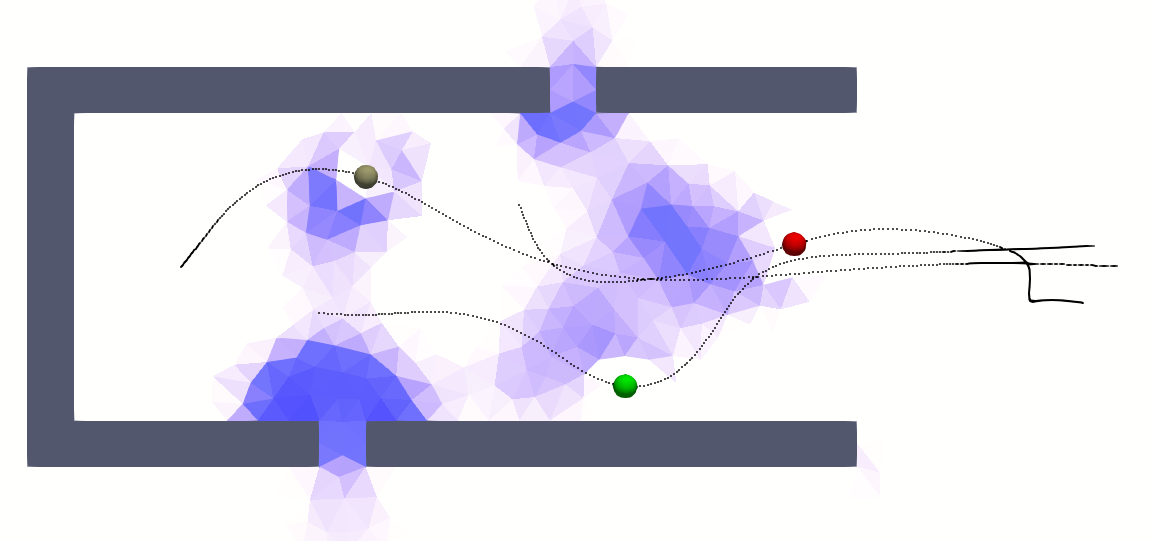}}
}
\subfloat[Solution at time $t=4.2$]{
\centering
\fbox{\includegraphics[width=.49\textwidth]{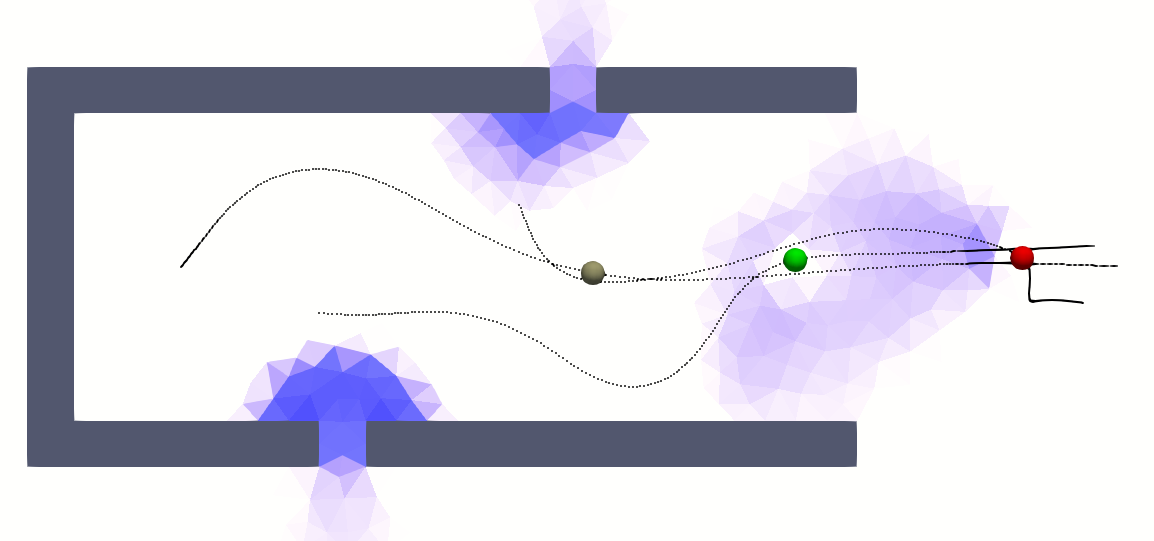}}
} \\
\subfloat[Solution at time $t=7.2$]{
\centering
\fbox{\includegraphics[width=.49\textwidth]{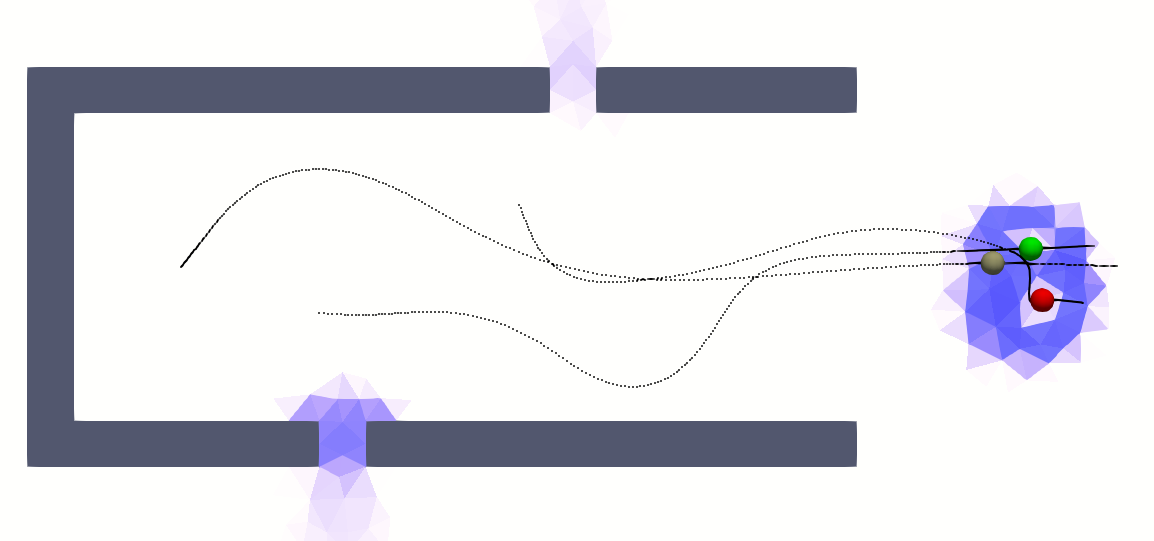}}
}
\subfloat[Agent intensities $c_i(t)$]{
\begin{tikzpicture}
\begin{axis}[enlargelimits=true, axis x line=center, axis y line=center,
		  axis line style = thick,
		  ytick distance=1,
		  width=0.48\textwidth, height=4.cm,
		  xlabel={\small $k$}, ylabel={\small $c_i$},
		  every axis x label/.style={
  		  	at={(ticklabel* cs:1.05)}, anchor=west,},
		  every axis y label/.style={
		    at={(ticklabel* cs:1.05)}, anchor=south,},
		  xmin=0, xmax=300, ymin=0, ymax=4.5, 
		  minor tick num=2, grid=both,
		  grid style={line width=.1pt, draw=gray!10}]
  	\addlegendimage{no markers, thick, red};
	\addlegendimage{no markers, thick, green};
	\addlegendimage{no markers, thick, color=agent3};	
	\addplot+[red, thick] table[col sep=space,
		 x expr=\coordindex, y=Ag_int_0,
		 mark=none]{experiment_2_agents/control.csv};
	\addlegendentry{$c_1(t)$};
	\addplot+[green, thick] table[col sep=space,
		 x expr=\coordindex, y=Ag_int_1,
		 mark=none]{experiment_2_agents/control.csv};
	\addlegendentry{$c_2(t)$};
	\addplot+[color=agent3, thick] table[col sep=space,
		 x expr=\coordindex, y=Ag_int_2,
		 mark=none]{experiment_2_agents/control.csv};
	\addlegendentry{$c_3(t)$};	
	\addplot+ [mark=none, thick, black, solid] coordinates {(0,  0) (0,  4.5)};
	\addplot+ [mark=none, thick, black, solid] coordinates {(40, 0) (40, 4.5)};
	\addplot+ [mark=none, thick, black, solid] coordinates {(80, 0) (80, 4.5)};
	\addplot+ [mark=none, thick, black, solid] coordinates {(140, 0) (140, 4.5)};
	\addplot+ [mark=none, thick, black, solid] coordinates {(240, 0) (240, 4.5)};
\end{axis}
\end{tikzpicture}}
\vspace{-2mm}
\begin{center}
\includegraphics[width=.4\textwidth]{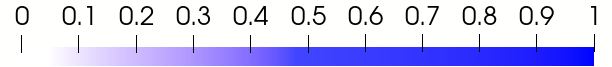}
\end{center}
\caption{Solution of the problem from \cref{sec:example_1}. The
colored background represents the density $\rho$; the dots are
the agent positions; the black curves are the agent trajectories.}
\label{fig:solution_example_1}
\end{figure}
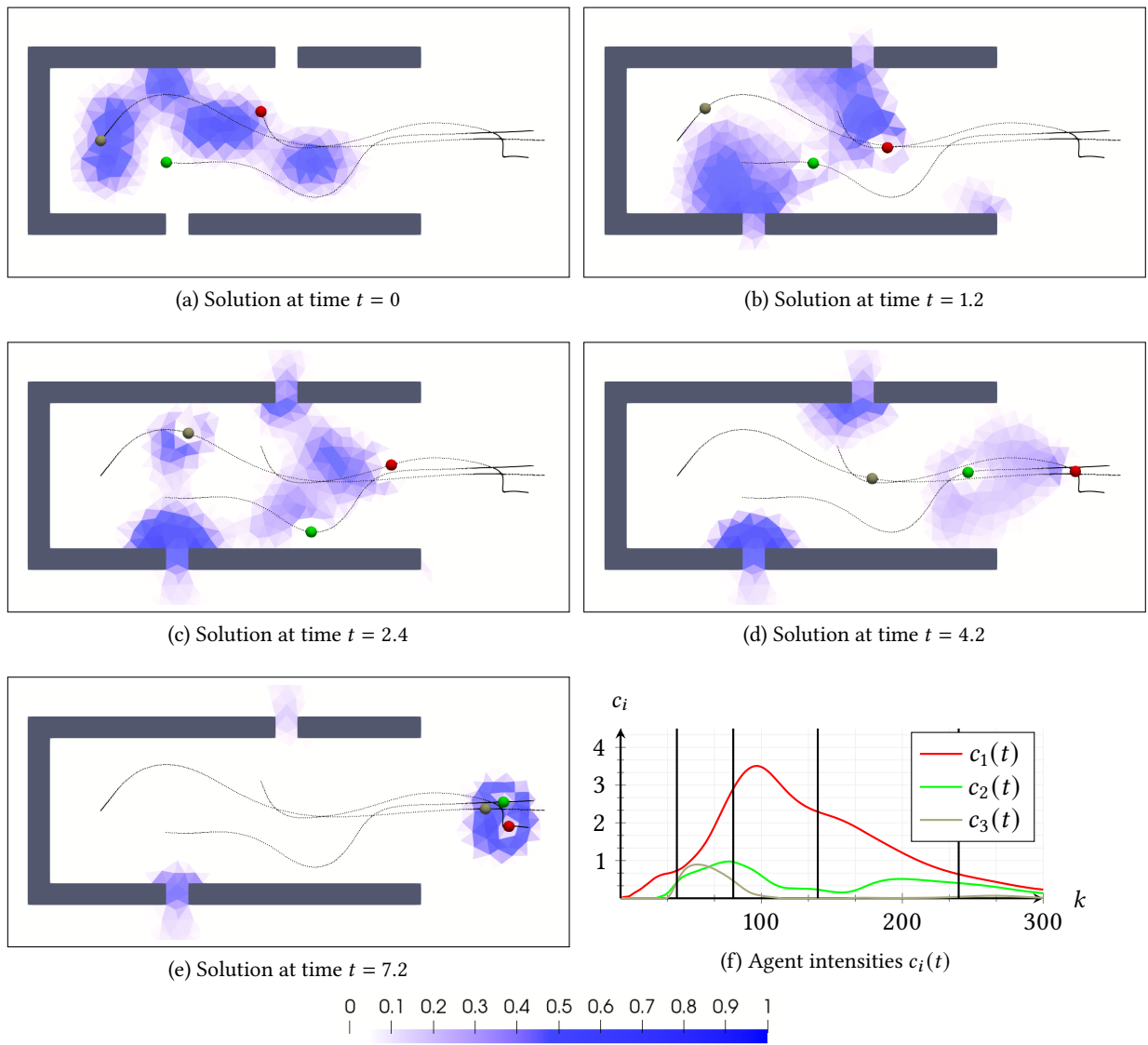

In a first numerical test we solve the problem
\cref{eq:optimal_control} in the domain $\Omega$ depicted in
\cref{fig:solution_example_1} with the following parameters.
\begin{align*}
	T&=9 & n_T &= 300 & \alpha_1 = \alpha_2 &= 5\cdot 10^{-2} & \gamma &= 10 & \zeta &= 10^{-2}
	& \mu &= 5\cdot 10^{-2}  \\
	\varepsilon &= 10^{-5} &
	\delta_1 &= 0.2 & \delta_2 &= 0.1 &
	\delta_3 &= 10^{-2} & \delta_4 &= 0.1.
\end{align*}
The initial density $\rho_0$ is the sum of 6 Gaussian bells, see also
\cref{fig:example_1_initial}.
The subdomain $\widetilde \Omega$ where densities are penalized is chosen to
cover the region within the walls.
Without any controlled agents, most of the people will squeeze
through the 2 smaller emergency exits in the south and north while
the large exit in the east is rarely used. To improve the evacuation 3
agents were introduced. The initial control $(\bu,\bc)$ was chosen in
such a way that the agent moves straight to the right outside of the
room having a constant the intensity.

\subsection{Example 2}
\label{sec:example_2}

\begin{figure}
\begin{center}
\subfloat[Solution at time $t=0$]{\fbox{\includegraphics[width=0.29\textwidth]{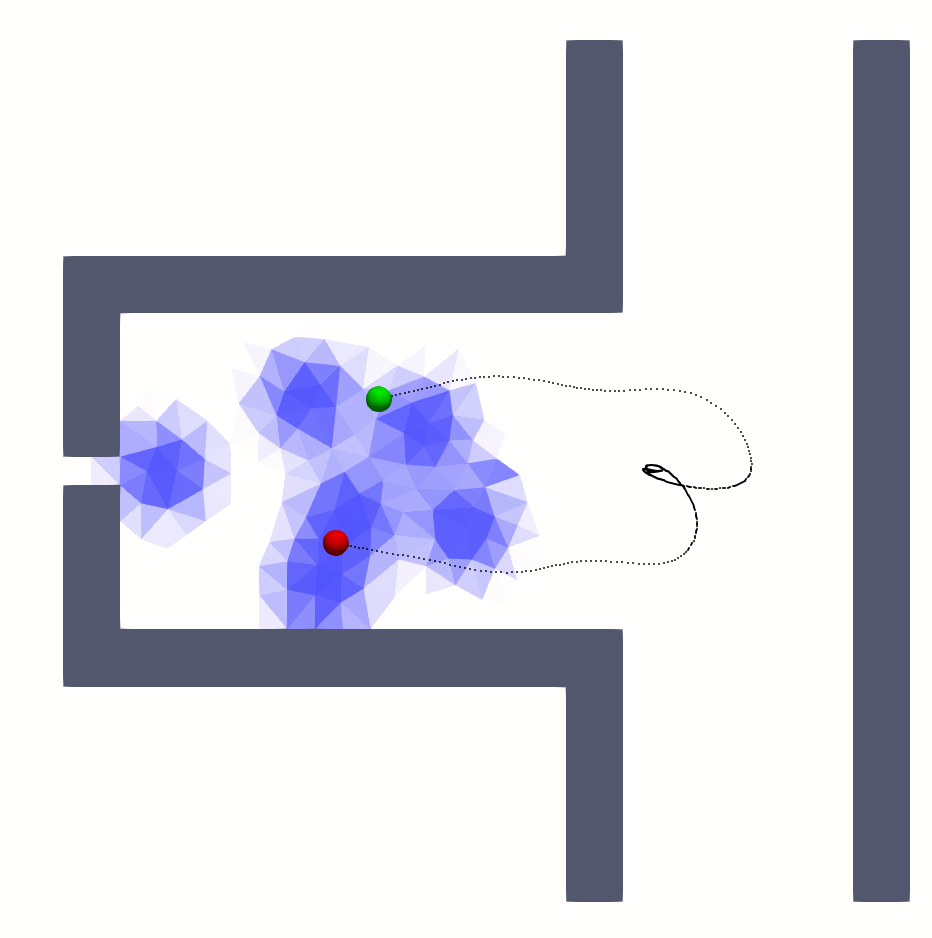}}}\quad
\subfloat[Solution at time $t=2.4$]{\fbox{\includegraphics[width=0.29\textwidth]{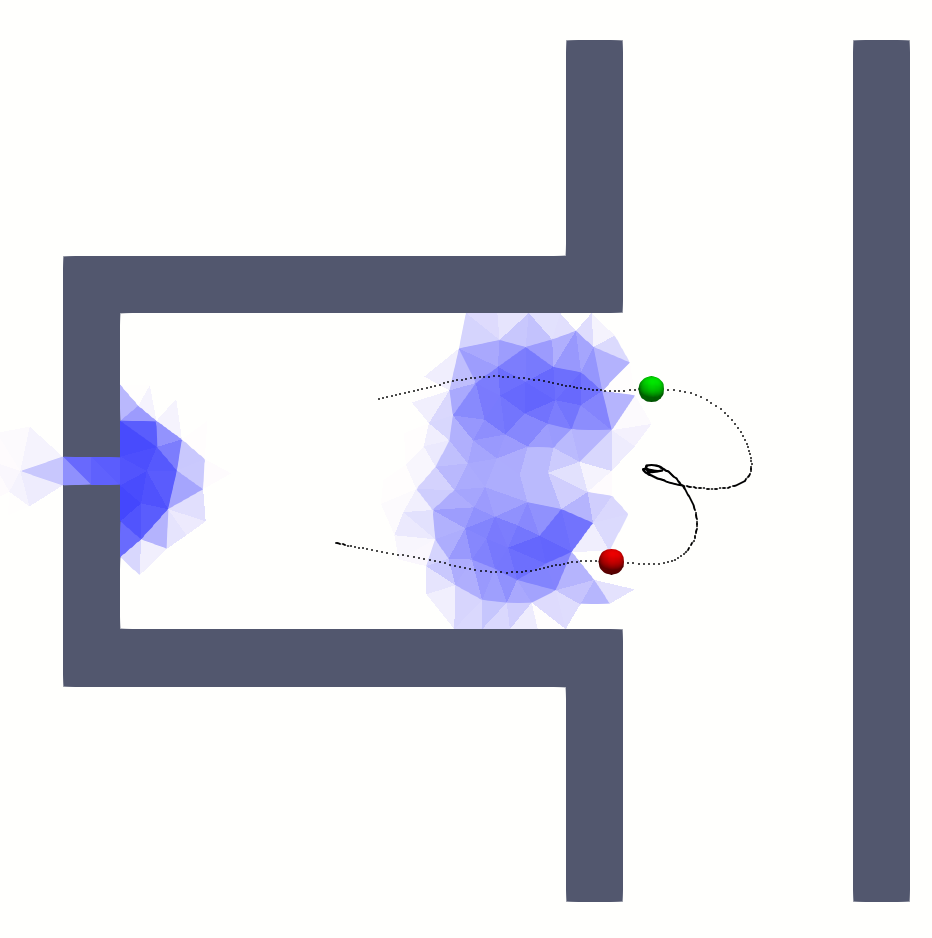}}}\quad
\subfloat[Solution at time $t=4.8$]{\fbox{\includegraphics[width=0.29\textwidth]{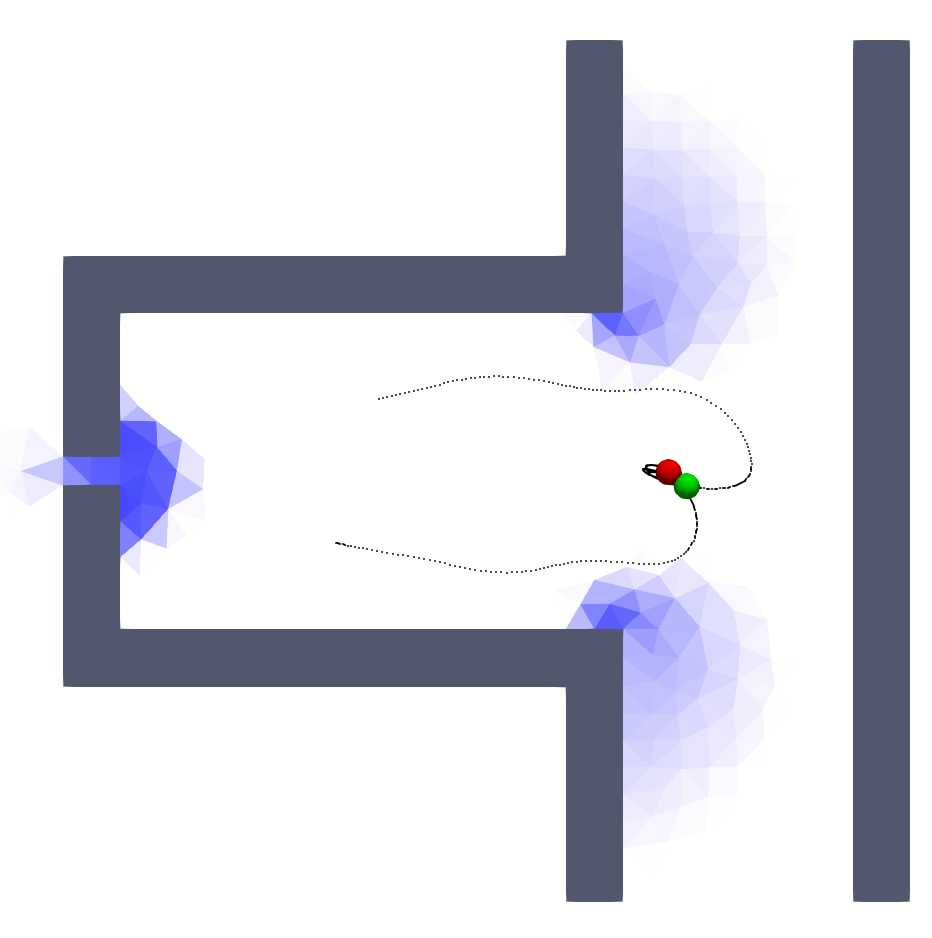}}}\\
\subfloat[Solution at time $t=8$]{\fbox{\includegraphics[width=0.29\textwidth]{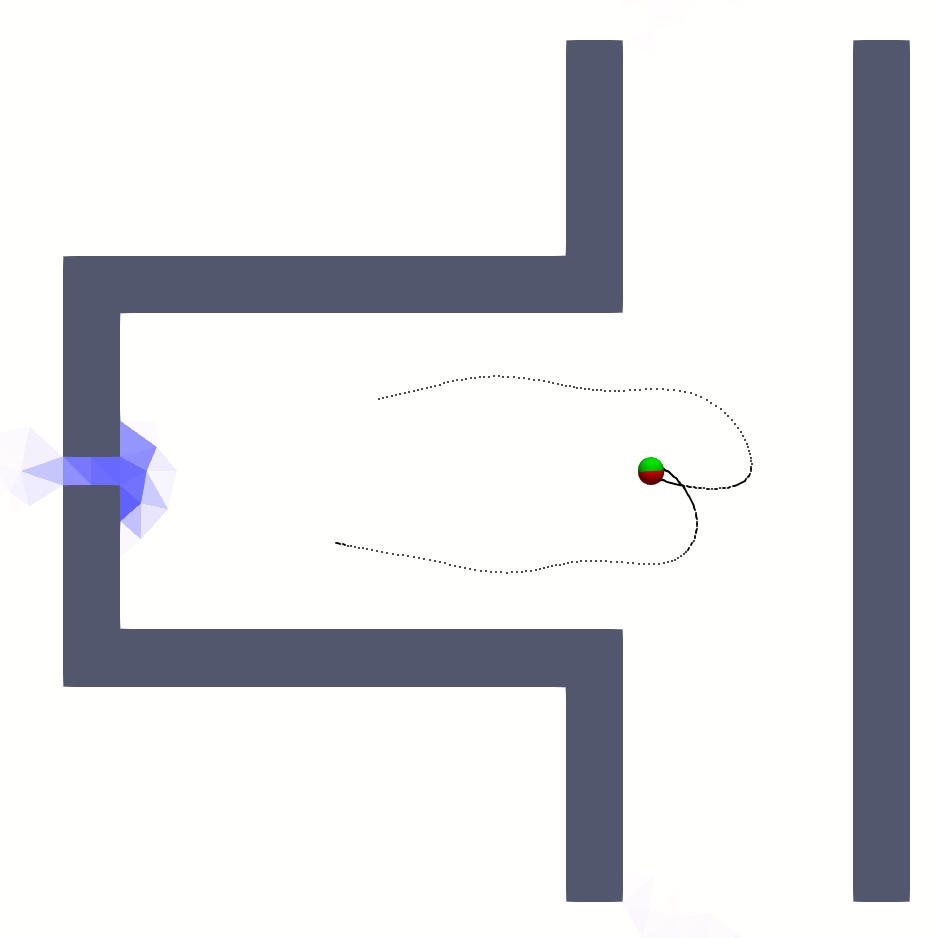}}}
\subfloat[Intensities $c_i$]{
\begin{tikzpicture}
\begin{axis}[enlargelimits=true, axis x line=center, axis y line=center,
		  axis line style = thick,
		  ytick distance=1,
		  width=0.65\textwidth, height=5cm,
		  every axis x label/.style={
  		  	at={(ticklabel* cs:1.05)}, anchor=west,},
		  every axis y label/.style={
		    at={(ticklabel* cs:1.05)}, anchor=south,},		  
		  xlabel={\small $k$}, ylabel={\small $c_i$},
		  xmin=0, xmax=220, ymin=0, ymax=6.5, 
		  minor tick num=2, grid=both,
		  grid style={line width=.1pt, draw=gray!10}]
  	\addlegendimage{no markers, thick, red};
	\addlegendimage{no markers, thick, yellow};
	\addplot+[red, thick] table[col sep=space,
		 x expr=\coordindex, y=Ag_int_0,
		 mark=none]{experiment_example_4/control.csv};
	\addlegendentry{$c_1(t)$};
	\addplot+[green, thick] table[col sep=space,
		 x expr=\coordindex, y=Ag_int_1,
		 mark=none]{experiment_example_4/control.csv};
	\addlegendentry{$c_2(t)$};
	\addplot+ [mark=none, thick, black, solid] coordinates {(0,  0) (0,  6.5)};
	\addplot+ [mark=none, thick, black, solid] coordinates {(60, 0) (60, 6.5)};
	\addplot+ [mark=none, thick, black, solid] coordinates {(120, 0) (120, 6.5)};
	\addplot+ [mark=none, thick, black, solid] coordinates {(200, 0) (200, 6.5)};	
\end{axis}
\end{tikzpicture}
}
\vspace{-2mm}
\begin{center}
\includegraphics[width=.4\textwidth]{./experiment_2_agents/color_legend.png}
\end{center}
\end{center}
\caption{Solution of the problem from \cref{sec:example_2} at various time steps $t_k$ and intensity of the agents.}
\label{fig:example2}
\end{figure}

In a second example we consider the domain illustrated in
\cref{fig:example2}. The initial density is concentrated near the slit
in the wall on the left-hand side. In an uncontrolled evacuation scenario the majority of
the people would leave the domain through this slit causing a massive
congestion leading to a very slow evacuation of the crowd.
The model and algorithm parameters chosen in the current example are as follows:
\begin{align*}
	T&=12 & n_T &= 300 & \alpha_1 = \alpha_2 &= 5\cdot 10^{-2} & \gamma &= 10 & \zeta &= 10^{-2}
	& \mu &= 5\cdot 10^{-2}  \\
	\varepsilon &= 10^{-5} &
	\delta_1 &= 0.2 & \delta_2 &= 0.1 &
	\delta_3 &= 10^{-2} & \delta_4 &= 0.1.
\end{align*}
This example shows that the evacuation can be significantly improved by using two agents with optimized trajectory and intensity. Interesting is, that the intensity is non-zero only in the time interval $t\in(0,3)$. The agents attract the people leading them sufficiently far away from the slit in the west and then they stop influencing the crowd. When being sufficiently far away from the slit the people find the way to the larger exits in the north and south on their own by using the movement direction determined by the potential $\phi$. 

\subsection{Example 3}
\label{sec:example_3}

In a third example we consider a square-shaped room with exits in the south, east and north. The exits have different width. The model and algorithm parameters are chosen as follows:
\begin{align*}
	T&=10 & n_T &= 300 & \alpha_1 = \alpha_2 &= 5\cdot 10^{-2} & \gamma &= 10 & \zeta &= 10^{-2}
	& \mu &= 5\cdot 10^{-2}  \\
	\varepsilon &= 10^{-5} &
	\delta_1 &= 0.2 & \delta_2 &= 0.1 &
	\delta_3 &= 10^{-2} & \delta_4 &= 0.1.
\end{align*}
The initial density is concentrated near the small exit and without a control of the crowd motion most of the people are blocking each other while squeezing through this small exit. Two agents were added in this scenario with the aim attracting the people in such a way that more of them find the other two exits in the north and east. The computed agent trajectories are quite short. It is interesting to observe that in the time interval $t\in (0,2)$ the agents just go to an optimal position sufficiently close to the crowd and attract them in the time interval $t\in (2,5)$, leading some of the people to the center of the room. At this point the agents drive their intensity to zero meaning that they stop influencing the crowd. However, when being sufficiently far away from the critical exit the people find the route to the less used exits on their own due to the movement rule determined by the potential $\phi$. 

\begin{figure}
\begin{center}
\subfloat[Solution at time $t=0$]{\fbox{\includegraphics[width=0.29\textwidth]{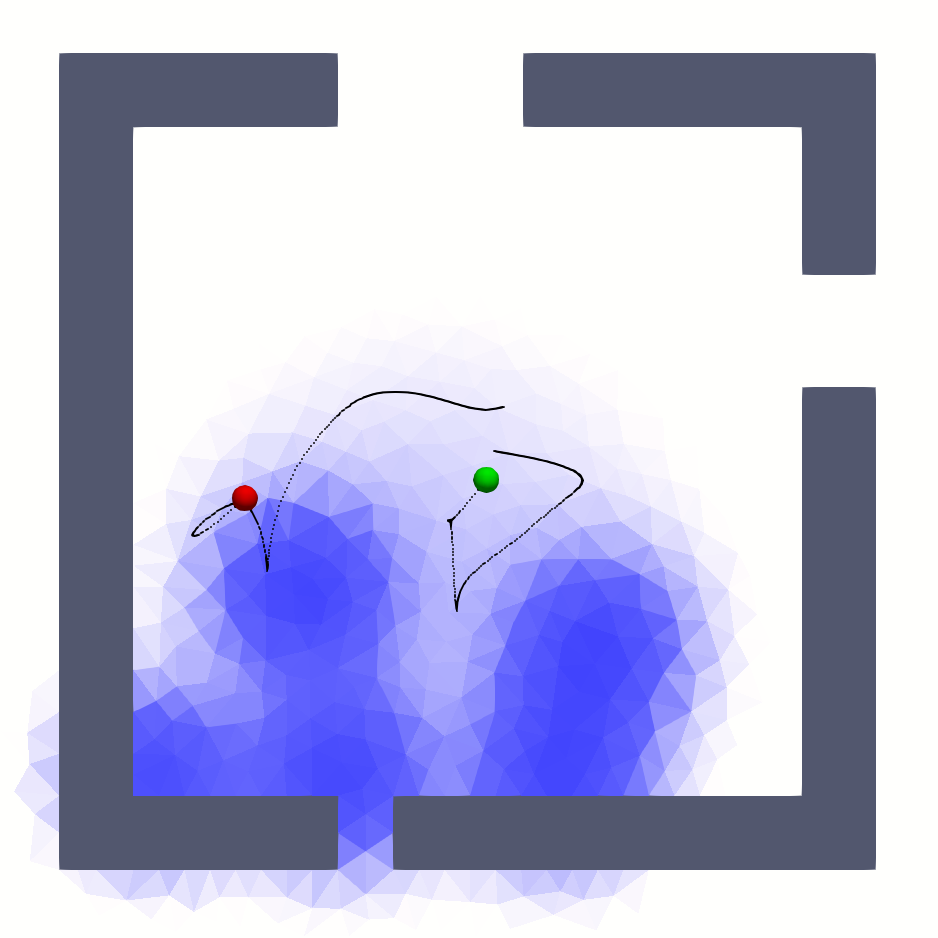}}}\quad
\subfloat[Solution at time $t=2$]{\fbox{\includegraphics[width=0.29\textwidth]{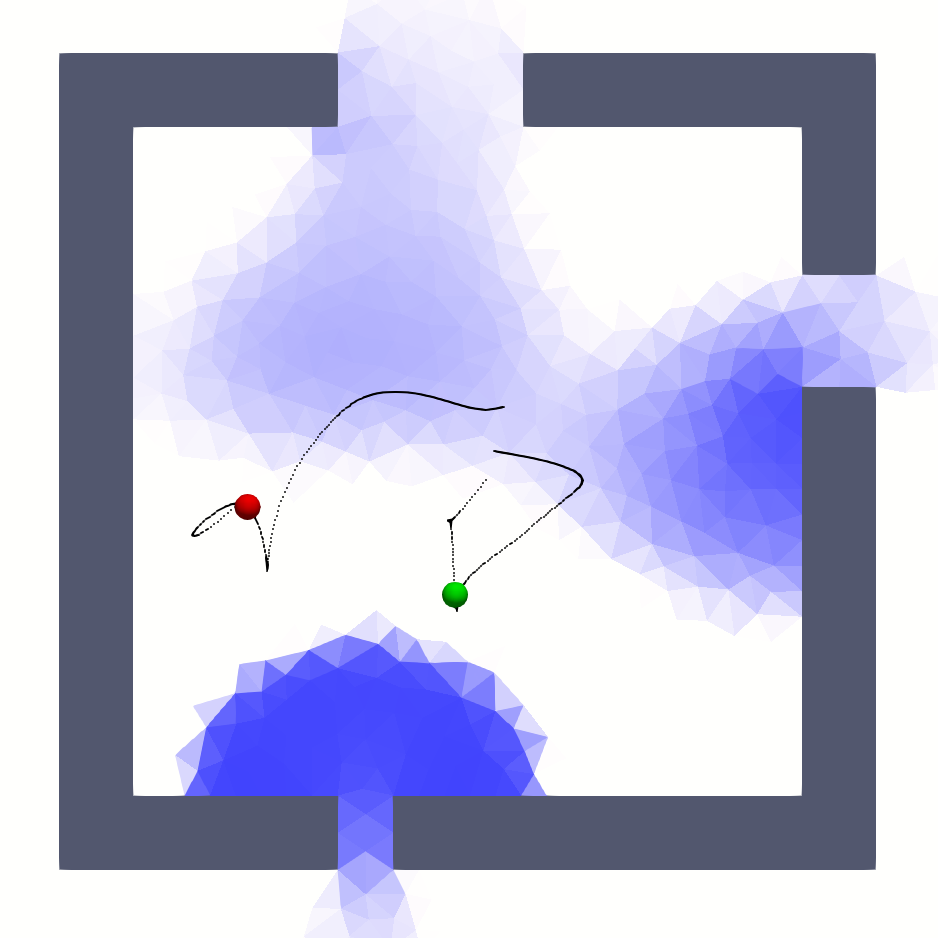}}}\quad
\subfloat[Solution at time $t=4$]{\fbox{\includegraphics[width=0.29\textwidth]{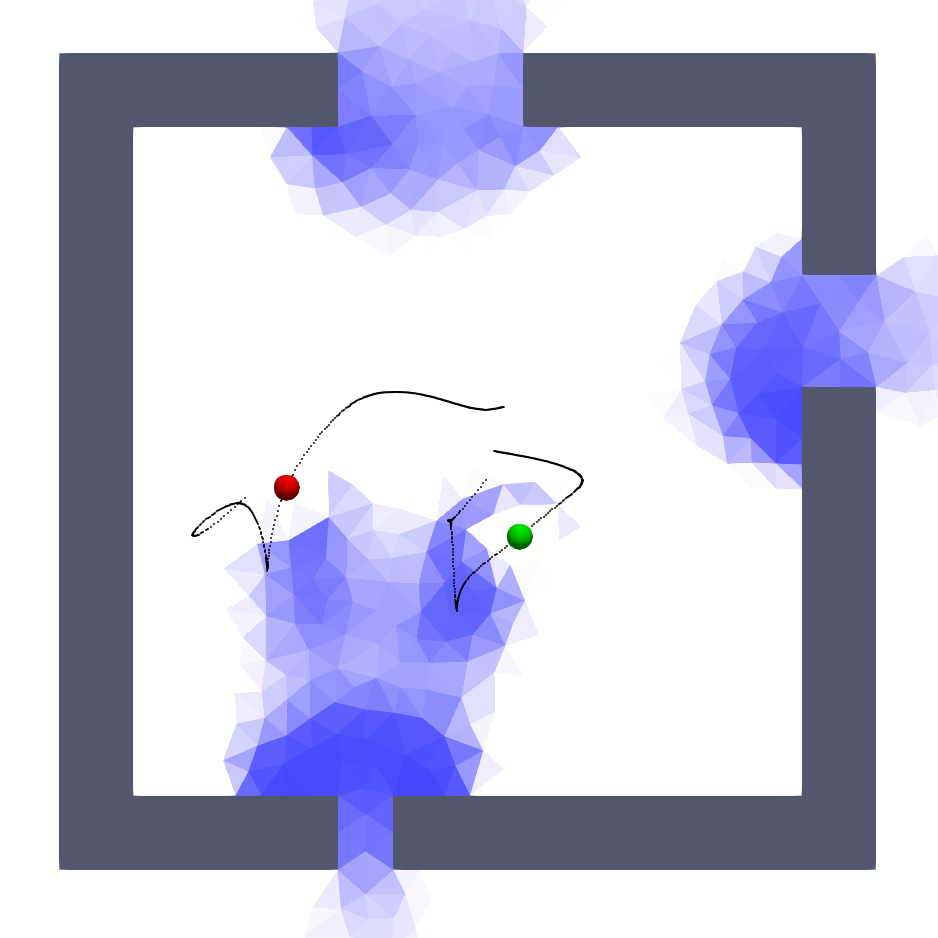}}}\\
\subfloat[Solution at time $t=6.7$]{\fbox{\includegraphics[width=0.29\textwidth]{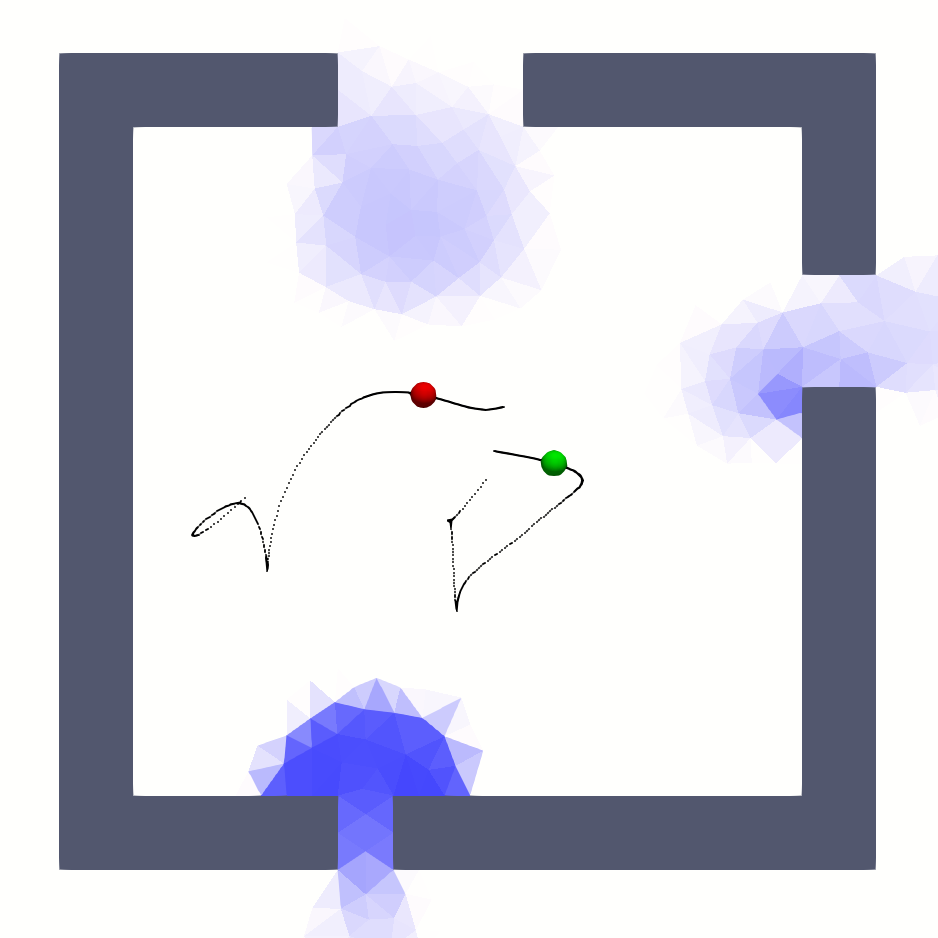}}}
\subfloat[Intensities $c_i$]{
\begin{tikzpicture}
\begin{axis}[enlargelimits=true, axis x line=center, axis y line=center,
		  axis line style = thick,
		  ytick distance=1,
		  width=0.65\textwidth, height=5cm,
		  xlabel={\small $k$}, ylabel={\small $c_i$},
		  every axis x label/.style={
  		  	at={(ticklabel* cs:1.05)}, anchor=west,},
		  every axis y label/.style={
		    at={(ticklabel* cs:1.05)}, anchor=south,},		  
		  xmin=0, xmax=300, ymin=0, ymax=2, 
		  minor tick num=2, grid=both,
		  grid style={line width=.1pt, draw=gray!10}]
  	\addlegendimage{no markers, thick, red};
	\addlegendimage{no markers, thick, green};
	\addplot+[red, thick] table[col sep=space,
		 x expr=\coordindex, y=Ag_int_0,
		 mark=none]{experiment_example_6/control.csv};
	\addlegendentry{$c_1(t)$};
	\addplot+[green, thick] table[col sep=space,
		 x expr=\coordindex, y=Ag_int_1,
		 mark=none]{experiment_example_6/control.csv};
	\addlegendentry{$c_2(t)$};
	\addplot+ [mark=none, thick, black, solid] coordinates {(0,  0) (0,  2)};
	\addplot+ [mark=none, thick, black, solid] coordinates {(60, 0) (60, 2)};
	\addplot+ [mark=none, thick, black, solid] coordinates {(120, 0) (120, 2)};
	\addplot+ [mark=none, thick, black, solid] coordinates {(200, 0) (200, 2)};	
\end{axis}
\end{tikzpicture}
}
\end{center}
\caption{Solution of the problem from \cref{sec:example_3} at various time steps $t_k$ and intensity of the agents.}
\label{fig:example3}
\end{figure}


\printbibliography

\end{document}